\documentclass[11pt]{article}
\usepackage{times,mathptm}
\usepackage{amsfonts,amscd,amssymb,theorem}
\input{diagrams}
\newtheorem{definition}{Definition}[section]
\newtheorem{lemma}[definition]{Lemma}
\newtheorem{proposition}[definition]{Proposition}
\newtheorem{corollary}[definition]{Corollary}
{\theorembodyfont{\rmfamily}\newtheorem{remark}[definition]{Remark}}
{\theorembodyfont{\rmfamily}\newtheorem{remarks}[definition]{Remarks}}
\newtheorem{theorem}[definition]{Theorem}
{\theorembodyfont{\rmfamily}}
{\theorembodyfont{\rmfamily}\newtheorem{examples}[definition]{Examples}}
{\theorembodyfont{\rmfamily}}

\newcommand{\End}{\rm{End}\,}
\newcommand{\Hom}{\rm{Hom}\,}

\def\rightact{\hbox{$\leftarrow$}}

\def\va{\varepsilon}
\def\v{\varphi}
\def\tl{\triangleleft}
\def\tr{\triangleright}
\def\rh{\rightharpoonup}
\def\lh{\leftharpoonup}

\def\ra{\rightarrow}
\def\a{\alpha}
\def\b{\beta}

\def\l{\lambda}
\def\r{\rho}
\def\cd{\cdot}

\def\O{\Omega}
\def\o{\omega}
\def\ov{\overline}
\def\un{\underline}
\def\mf{\mathfrak}
\def\mb{\mathbb}

\newcommand{\mfa}{\mbox{$\mf {a}$}}
\newcommand{\mfb}{\mbox{$\mf {b}$}}

\newcommand{\hba}{\mbox{$_H{\cal M}_{\mf A}^H$}}

\newcommand{\smi}{\mbox{$S^{-1}$}}
\newcommand{\gsm}{\mbox{$\blacktriangleright \hspace{-0.7mm}<$}}
\newcommand{\gtl}{\mbox{${\;}$$>\hspace{-1.5mm}\triangleleft$${\;}$}}
\newcommand{\trl}{\mbox{${\;}$$\triangleright \hspace{-1.5mm}<$${\;}$}}
\newcommand{\ovsm}{\mbox{${\;}$$\ov {\#}$${\;}$}}
\newcommand{\tx}{\mbox{$\tilde {x}$}}
\newcommand{\tX}{\mbox{$\tilde {X}$}}
\newcommand{\ty}{\mbox{$\tilde {y}$}}
\newcommand{\tY}{\mbox{$\tilde {Y}$}}
\newcommand{\tz}{\mbox{$\tilde {z}$}}

\newcommand{\tpra}{\mbox{$\tilde {p}^1_{\rho }$}}
\newcommand{\tprb}{\mbox{$\tilde {p}^2_{\rho }$}}
\newcommand{\tPra}{\mbox{$\tilde {P}^1_{\rho }$}}
\newcommand{\tPrb}{\mbox{$\tilde {P}^2_{\rho }$}}
\newcommand{\tqra}{\mbox{$\tilde {q}^1_{\rho }$}}
\newcommand{\tqrb}{\mbox{$\tilde {q}^2_{\rho }$}}
\newcommand{\tQra}{\mbox{$\tilde {Q}^1_{\rho }$}}
\newcommand{\tQrb}{\mbox{$\tilde {Q}^2_{\rho }$}}
\newcommand{\und}{\mbox{$\un {\Delta }$}}
\newcommand{\une}{\mbox{$\un {\va }$}}
\newcommand{\una}{\mbox{$c_{\un {1}}$}}
\newcommand{\unb}{\mbox{$c_{\un {2}}$}}
\newcommand{\unaa}{\mbox{$c_{(\un {1}, \un {1})}$}}
\newcommand{\unab}{\mbox{$c_{(\un {1}, \un {2})}$}}

\newcommand{\ba}{\mbox{$(\mb {A}, \l , \r , \Phi _{\l },
\Phi _{\r }, \Phi _{\l , \r })$}}
\newcommand{\tstc}{\mbox{$_H^C{\cal M}_{\mb {A}}^H$}}
\newcommand{\mbA}{\mbox{$\mb {A}$}}
\def\rawo\lonra{\longrightarrow}

\def\ot{\otimes}

\newcommand{\thlabel}[1]{\label{th:#1}}
\newcommand{\thref}[1]{Theorem~\ref{th:#1}}
\newcommand{\selabel}[1]{\label{se:#1}}
\newcommand{\seref}[1]{Section~\ref{se:#1}}
\newcommand{\lelabel}[1]{\label{le:#1}}
\newcommand{\leref}[1]{Lemma~\ref{le:#1}}
\newcommand{\prlabel}[1]{\label{pr:#1}}
\newcommand{\prref}[1]{Proposition~\ref{pr:#1}}
\newcommand{\colabel}[1]{\label{co:#1}}
\newcommand{\coref}[1]{Corollary~\ref{co:#1}}
\newcommand{\relabel}[1]{\label{re:#1}}
\newcommand{\reref}[1]{Remark~\ref{re:#1}}
\newcommand{\exlabel}[1]{\label{ex:#1}}

\newcommand{\delabel}[1]{\label{de:#1}}

\newcommand{\eqlabel}[1]{\label{eq:#1}}
\newcommand{\eqref}[1]{(\ref{eq:#1})}

\newenvironment{proof}{{\it Proof.}}{\hfill $ \square $ \vskip 4mm}
\begin{document}
\title{Two-sided (two-cosided) Hopf modules and
Doi-Hopf modules for quasi-Hopf algebras\thanks{Research supported by the
bilateral
project ``Hopf Algebras in Algebra, Topology, Geometry and Physics" of the
Flemish and Romanian governments.}}
\author{D. Bulacu\thanks{This paper was written while the first author was
visiting the Free University of Brussels, VUB (Belgium); he would
like to thank VUB for its warm hospitality.}\\ Faculty of
Mathematics\\ University of Bucharest\\ RO-70109 Bucharest 1,
Romania\and S. Caenepeel
\\ Faculty of Applied Sciences\\
Free University of Brussels, VUB\\ B-1050 Brussels, Belgium}
\date{}
\maketitle

\begin{abstract}
Let $H$ be a finite dimensional quasi-Hopf algebra over a field
$k$ and ${\mf A}$ a right $H$-comodule algebra in the sense of
\cite{hn1}. We first show that on the $k$-vector space ${\mf A}\ot
H^*$ we can define an algebra structure, denoted by ${\mf A}\ovsm
H^*$, in the monoidal category of left $H$-modules (i.e. ${\mf
A}\ovsm H^*$ is an $H$-module algebra in the sense of \cite{bpv}).
Then we will prove that the category of two-sided $({\mf A},
H)$-bimodules $\hba $ is isomorphic to the category of relative
$({\mf A}\ovsm H^*, H^*)$-Hopf modules, as introduced in
\cite{bn}. In the particular case where ${\mf A}=H$, we will
obtain the Nill's result announced in \cite{hn3}. We will also
introduce the categories of Doi-Hopf modules and two-sided
two-cosided Hopf modules and we will show that they are in certain
situations isomorphic to module categories.
\end{abstract}
%%%%%%%%%%%%%%%%%%%%%%%%%%%%%%%%%%%%%%%%%%%%%%%%%%%
\section{Introduction}\selabel{0}
%%%%%%%%%%%%%%%%%%%%%%%%%%%%%%%%%%%%%%%%%%%%%%%%%%%%
Quasi-bialgebras and quasi-Hopf algebras were introduced by
Drinfeld \cite{d1} in connection with the Knizhnik-Zamolodchikov
equations \cite{k}. Let $k$ be a field, $H$ an associative algebra
and $\Delta :\ H\ra H\ot H$ and $\va :\ H\ra k$ two algebra
morphisms. Roughly speaking, $H$ is a quasi-bialgebra if  the
category $_H{\cal M}$ of left $H$-modules, equipped with the
tensor product of vector spaces endowed with the diagonal
$H$-module structure given via $\Delta $, and with unit object $k$
viewed as a left $H$-module via $\va $, is a monoidal category.
The comultiplication $\Delta$ is not coassociative but is
quasi-coassociative in the sense that $\Delta $ is coassociative
up to conjugation by an invertible element $\Phi \in H\ot H\ot H$.
Moreover, $H$ is a quasi-Hopf algebra if and only if each finite
dimensional left $H$-module has a dual $H$-module. Note that, the
definition of a quasi-bialgebra and a quasi-Hopf
algebra is not self dual.\\%
Since $H$ is not a coassociative coalgebra, it is impossible to
define comodules over $H$. However, since $H$ is a coalgebra in
the monoidal category of $(H,H)$-bimodules $_H{\cal M}_H$, we can
define Hopf modules in $_H{\cal M}_H^H$. More exactly, a Hopf
module in $_H{\cal M}_H^H$ is a right $H$-comodule in the monoidal
category $_H{\cal M}_H$, cf. \cite{hn3}. Also, we can define a
Hopf module in $_H^H{\cal M}_H^H$: a Hopf module in $_H^H{\cal
M}_H^H$ is an $H$-$H$-bicomodule in the monoidal category $_H{\cal
M}_H$, cf. \cite{sh}. \\%
Using the theory of (co)algebras and (co)modules in monoidal
categories we can define categories of relative Hopf modules. If
$H$ is a finite dimensional quasi-bialgebra and $A$ an algebra in
the monoidal category $_H{\cal M}$, then a relative Hopf module in
${\cal M}_A^{H^*}$ is a right $H^*$-comodule ($H^*$ is a
coassociative coalgebra) which is also a right $A$-module in the
monoidal category ${\cal M}^{H^*}$, see \cite{bn}. Similarly, if
we start with a coalgebra $C$ in the monoidal category of right
$H$-modules ${\cal M}_H$, then a relative Hopf module in
${}^C{\cal M}_H$ is a left $C$-comodule in the monoidal category
${\cal M}_H$. Of course, when $H$ is an ordinary bialgebra these
categories are exactly the ones defined by Doi in \cite{doi2}.
Unfortunately, at the moment, the basic category of Hopf modules
${\cal M}_H^H$ is not defined (a possible way to define this
category will be presented later in Section 3).\\%
In this paper, our goal is to generalize all the categories
presented above and to prove that, in the finite dimensional case,
they are isomorphic to module categories. Essential tools will be
the notions of comodule algebra and bicomodule algebra over a
quasi-bialgebra $H$ \cite{hn1}. An associative algebra ${\mf A}$
is called a right $H$-comodule algebra if there exists an algebra
map $\r : {\mf A}\ra {\mf A}\ot H$ which is ``almost" associative,
in a way similar to the ``almost" coassociativity of the
comultiplication on a quasi-bialgebra; a detailed definition will
be presented below. In a similarl way, we can define the notion of
left $H$-comodule algebra. If $H$ is a finite dimensional
quasi-bialgebra and ${\mf A}$ a right $H$-comodule algebra, then
we can define a multiplication on the $k$-vector space ${\mf A}\ot
H^*$; we obtain an algebra in the monoidal category $_H{\cal M}$
(i.e. a left $H$-module algebra), which will be denoted by ${\mf
A}\ovsm H^*$. Notice that, in the Hopf case, ${\mf A}\ovsm H^*$ is
just the usual smash product of ${\mf A}$ and $H^*$, and this is
why we call the monoidal algebra ${\mf A}\ovsm H^*$ the
quasi-smash product of ${\mf A}$ and $H^*$. Also remark that, in
the quasi-Hopf case, ${\mf A}\ovsm H^*$ is not an associative
algebra because the associativity constraints of $_H{\cal M}$ are
not the trivial ones but the ones given by the reassociator $\Phi
$ of $H$.\\ In \seref{3} we will define the category of two-sided
$({\mf A}, H)$-Hopf modules $_H{\cal M}_{{\mf A}}^H$. Concerning
its description in terms of relative Hopf modules, a central role
is played by ${\mf A}\ovsm H^*$. First, the definition of a
two-sided $({\mf A}, H)$-Hopf module is slightly more general than
the one given for a Hopf module in $_H{\cal M}_H^H$. Its objects
are $H$-${\mf A}$-bimodules and ``almost" right $H$-comodules such
that, in the usual way, the right $H$-coaction is left $H$-linear
an a right ${\mf A}$-linear map. If ${\mf A}=H$, then $_H{\cal
M}_{{\mf A}}^H$ is just the category of Hopf modules ${}_H{\cal
M}_H^H$ described above. Secondly, if $H$ is a finite dimensional
quasi-Hopf algebra then $_H{\cal M}_{{\mf A}}^H$ is isomorphic to
the category of relative $({\mf A}\ovsm H^*, H^*)$-Hopf modules
(this generalizes \cite[Proposition 2.3]{drvo}). Following
\cite{bn}, we can show in this case that ${}_H{\cal M}_{{\mf
A}}^H$ is isomorphic to the category of right $({\mf A}\ovsm
H^*)\# H$-modules, where $({\mf A}\ovsm H^*)\# H$ denotes the
smash product algebra (in the sense of \cite{bpv}) of ${\mf
A}\ovsm H^*$ and $H$. In particular, taking ${\mf A}=H$ we recover
Nill's result, as announced in \cite{hn3}, which states that
${}_H{\cal M}_H^H$ is isomorphic to the category of right modules
over the two-sided crossed product $H\gtl H^*\trl H$. In \seref{4}
we will prove that the two-sided crossed product constructed in
\cite{hn1} is in fact a generalized smash product. As a
consequence, $(H\ovsm H^*)\# H$  is just the two-sided crossed
product $H\gtl H^*\trl H$ (as an algebra).\\%
In the second part of the paper we will study the category of
two-sided two-cosided Hopf modules $\tstc$. Here $C$ is a
coalgebra in the monoidal category of $(H,H)$-bimodules ${}_H{\cal
M}_H$ (i.e. an $H$-bimodule coalgebra), and ${\mb A}$ is an
$H$-bicomodule algebra in the sense of \cite{hn1}. Roughly
speaking, an object in $\tstc $ is a two-sided $({\mb A}, H)$-Hopf
module which is also an ``almost" left $C$-comodule such that the
left $C$-coaction is compatible with the other structure maps. In
\seref{5} we will show that if $C$ and $H$ are finite dimensional
then $\tstc $ is isomorphic to a category of right modules. To
this end we will describe first $\tstc $ as a category of Doi-Hopf
modules. If ${\mf B}$ is a left $H$-comodule algebra and $C$ is a
right $H$-module coalgebra then the category of right-left $(H,
{\mf B}, C)$-Doi-Hopf modules $^C{\cal M}(H)_{\mf B}$ is a
straightforward generalization of the category of relative Hopf
modules $^C{\cal M}_H$. When $C$ is finite dimensional, $^C{\cal
M}(H)_{\mf B}$ is isomorphic to the category of right modules over
the generalized smash product $C^*\gsm {\mf B}$. Now, returning to
the category $\tstc $, if $H$ is finite dimensional then we will
show that $(\mbA \ovsm H^*)\# H$ is a left $H\ot H^{\rm
op}$-comodule algebra (here "op" means the opposite multiplication
on $H$) so, it makes sense to consider the category of Doi-Hopf
modules $^C{\cal M}(H\ot H^{\rm op})_{(\mbA \ovsm H^*)\# H}$.
Moreover, the main results asserts that $\tstc $ is isomorphic to
$^C{\cal M}(H\ot H^{\rm op})_{(\mbA \ovsm H^*)\# H}$ (this
generalizes \cite[Proposition 2.3]{bdr}). In particular, if $C$ is
also finite dimensional, we obtain that $\tstc $ is isomorphic to
the category of right modules over the generalized smash product
${\cal A}=C^*\gsm ((\mbA \ovsm H^*)\# H)$. In the Hopf case, the
left-handed version of this result was first obtained by Cibils
and Rosso \cite{cr}. More precisely, they define an algebra $X$
having the property that the category ${}_{H^*}^{H^*}{\cal
M}_{H^*}^{H^*}$ is isomorphic to the category of left $X$-modules.
Recently, Panaite \cite{p} introduced two other algebras $Y$ and
$Z$ with the same property as $X$. $Y$ is the two-sided crossed
product $H^*\# (H\ot H^{\rm op})\# H^{*\rm op}$ and $Z$ is the
diagonal crossed product (in the sense of \cite{hn1}) $(H^*\ot
H^{*\rm op})\Join (H\ot H^{\rm op})$.

%%%%%%%%%%%%%%%%%%%%%%%%%%%%%%%%%%%%%%%%%%%%%%%%%%%%%%%%%%%%%%
\section{Preliminary results}\selabel{1}
%%%%%%%%%%%%%%%%%%%%%%%%%%%%%%%%%%%%%%%%%%%%%%%%%%%%%%%%%%%%%%%
\subsubsection*{Quasi-Hopf algebras}
\setcounter{equation}{0}
We work over a commutative field $k$. All algebras, linear spaces
etc. will be over $k$; unadorned $\ot $ means $\ot_k$. Following
Drinfeld \cite{d1}, a quasi-bialgebra is a fourtuple $(H, \Delta ,
\va , \Phi )$ where $H$ is an associative algebra with unit,
$\Phi$ is an invertible element in $H\ot H\ot H$, and $\Delta :\
H\ra H\ot H$ and $\va :\ H\ra k$ are algebra homomorphisms
satisfying the identities
\begin{eqnarray}
&&(id \ot \Delta )(\Delta (h))=%
\Phi (\Delta \ot id)(\Delta (h))\Phi ^{-1},\label{q1}\\[1mm]%
&&(id \ot \va )(\Delta (h))=h\ot 1, %
\mbox{${\;\;\;}$}%
(\va \ot id)(\Delta (h))=1\ot h,\label{q2}
\end{eqnarray}
for all $h\in H$, and $\Phi$ has to be a normalized $3$-cocycle,
in the sense that
\begin{eqnarray}
&&(1\ot \Phi)(id\ot \Delta \ot id) (\Phi)(\Phi \ot 1)= (id\ot id
\ot \Delta )(\Phi ) (\Delta \ot id \ot id)(\Phi
),\label{q3}\\[1mm]%
&&(id \ot \va \ot id )(\Phi )=1\ot 1\ot 1.\label{q4}
\end{eqnarray} The map $\Delta $ is called the coproduct or the
comultiplication, $\va $ the counit and $\Phi $ the reassociator.
As for Hopf algebras \cite{sw} we denote $\Delta (h)=\sum h_1\ot
h_2$, but since $\Delta$ is only quasi-coassociative we adopt the
further
convention%
$$
(\Delta \ot id)(\Delta (h))= \sum h_{(1, 1)}\ot h_{(1, 2)}\ot h_2,
\mbox{${\;\;\;}$} (id\ot \Delta )(\Delta (h))=
\sum h_1\ot h_{(2, 1)}\ot h_{(2,2)}, %
$$
for all $h\in H$. We will denote the tensor components of $\Phi$
by capital letters, and the ones of $\Phi^{-1}$ by small letters,
namely
\begin{eqnarray*}
&&\Phi=\sum X^1\ot X^2\ot X^3= \sum T^1\ot T^2\ot T^3= \sum V^1\ot
V^2\ot V^3=\cdots\\
&&\Phi^{-1}=\sum x^1\ot x^2\ot x^3= \sum
t^1\ot t^2\ot t^3= \sum v^1\ot v^2\ot v^3=\cdots
\end{eqnarray*}
$H$ is called a quasi-Hopf algebra if, moreover, there exists an
anti-automorphism $S$ of the algebra $H$ and elements $\a , \b \in
H$ such that, for all $h\in H$, we have:
\begin{eqnarray}
&&\sum S(h_1)\a h_2=\va (h)\a \mbox{${\;\;\;}$ and ${\;\;\;}$}
\sum h_1\b S(h_2)=\va (h)\b ,\label{q5}\\[1mm]%
&&\sum X^1\b S(X^2)\a X^3=1 %
\mbox{${\;\;\;}$ and${\;\;\;}$}%
\sum S(x^1)\a x^2\b S(x^3)=1.\label{q6}
\end{eqnarray}
For a quasi-Hopf algebra the antipode is determined
uniquely up to a transformation $\a \mapsto U\a $, $\b \mapsto \b
U^{-1}$, $S(h)\mapsto US(h)U^{-1}$, where $U\in H$ is invertible.
The axioms for a quasi-Hopf algebra imply that $\va (\a )\va (\b
)=1$, so, by rescaling $\a $ and $\b $, we may assume without loss
of generality that $\va (\a )=\va (\b )=1$ and $\va \circ S=\va $.
The identities (\ref{q2}), (\ref{q3}) and (\ref{q4}) also imply
that
\begin{equation}\label{q7}
(\va \ot id\ot id)(\Phi )= (id \ot id\ot \va )(\Phi )=1\ot 1\ot 1.
\end{equation}
Next we recall that the definition of a quasi-Hopf algebra is
``twist coinvariant" in the following sense. An invertible element
$F\in H\ot H$ is called a {\sl gauge transformation} or {\sl
twist} if $(\va \ot id)(F)=(id\ot \va)(F)=1$. If $H$ is a
quasi-Hopf algebra and $F=\sum F^1\ot F^2\in H\ot H$ is a gauge
transformation with inverse $F^{-1}=\sum G^1\ot G^2$, then we can
define a new quasi-Hopf algebra $H_F$ by keeping the
multiplication, unit, counit and antipode of $H$ and replacing the
comultiplication, antipode and the elements $\alpha$ and $\beta$
by
\begin{eqnarray}
&&\Delta _F(h)=F\Delta (h)F^{-1},\label{g1}\\[1mm]%
&&\Phi_F=(1\ot F)(id \ot \Delta )(F) \Phi (\Delta \ot id)
(F^{-1})(F^{-1}\ot 1),\label{g2}\\[1mm]%
&&\a_F=\sum S(G^1)\a G^2,%
\mbox{${\;\;\;}$}%
\b_F=\sum F^1\b S(F^2).\label{g3}
\end{eqnarray}
It is well-known that the antipode of a Hopf
algebra is an anti-coalgebra morphism. For a quasi-Hopf algebra,
we have the following statement: there exists a gauge
transformation $f\in H\ot H$ such that
\begin{equation} \label{ca}
f\Delta (S(h))f^{-1}= \sum (S\ot S)(\Delta ^{\rm op}(h))
\mbox{,${\;\;\;}$for all $h\in H$,}
\end{equation}
where $\Delta ^{\rm op}(h)=\sum h_2\ot h_1$. $f$ can be computed
explicitly. First set
\begin{equation}
\sum A^1\ot A^2\ot A^3\ot A^4= (\Phi \ot 1) (\Delta \ot id\ot
id)(\Phi ^{-1}),
\end{equation}
\begin{equation} \sum B^1\ot B^2\ot B^3\ot B^4=
(\Delta \ot id\ot id)(\Phi )(\Phi ^{-1}\ot 1)
\end{equation}
and then define $\gamma, \delta\in H\ot H$ by
\begin{equation} \label{gd}%
\gamma =\sum S(A^2)\a A^3\ot S(A^1)\a A^4~~{\rm and}~~ \delta
=\sum B^1\b S(B^4)\ot B^2\b S(B^3).
\end{equation}
$f$ and $f^{-1}$ are then given by the formulas
\begin{eqnarray}
f&=&\sum (S\ot S)(\Delta ^{\rm op}(x^1)) \gamma \Delta (x^2\b
S(x^3)),\label{f}\\
f^{-1}&=&\sum \Delta (S(x^1)\a x^2) \delta (S\ot S)(\Delta
^{\rm op}(x^3)).\label{g}
\end{eqnarray}
$f$ satisfies the following relations:
\begin{equation} \label{gdf}%
f\Delta (\a )=\gamma , \mbox{${\;\;\;}$} %
\Delta (\b )f^{-1}=\delta .
\end{equation}
Furthermore the corresponding twisted reassociator (see
(\ref{g2})) is given by
\begin{equation} \label{pf}
\Phi _f=\sum (S\ot S\ot S)(X^3\ot X^2\ot X^1).
\end{equation}
In a Hopf algebra $H$, we obviously have the identity $$\sum
h_1\ot h_2S(h_3)=h\ot 1,~{\rm for~all~}h\in H.$$ We will need the
generalization of this formula to the quasi-Hopf algebra setting.
Following \cite{hn1}, \cite{hn2}, we define
\begin{equation} \label{qr}
p_R=\sum p^1_R\ot p^2_R=\sum x^1\ot x^2\b S(x^3),%
\mbox{${\;\;\;}$}%
q_R=\sum q^1_R\ot q^2_R=\sum X^1\ot S^{-1}(\a X^3)X^2,%
\end{equation}
\begin{equation}\label{ql}
p_L=\sum p^1_L\ot p^2_L=\sum X^2\smi (X^1\b )\ot X^3,%
\mbox{${\;\;\;}$}%
q_L=\sum q^1_L\ot q^2_L=\sum S(x^1)\a x^2\ot x^3.%
\end{equation}
For all $h\in H$, we then have
\begin{equation} \label{qr1}
\sum \Delta (h_1)p_R[1\ot S(h_2)]=p_R[h\ot 1],
\mbox{${\;\;\;}$}%
\sum [1\ot S^{-1}(h_2)]q_R\Delta (h_1)=(h\ot 1)q_R,%
\end{equation}
\begin{equation}\label{ql1}
\sum \Delta (h_2)p_L[\smi (h_1)\ot 1]=p_L(1\ot h),%
\mbox{${\;\;\;}$}%
\sum [S(h_1)\ot 1]q_L\Delta (h_2)=(1\ot h)q_L,
\end{equation}
and
\begin{equation} \label{pqr}
\sum \Delta (q^1_R)p_R[1\ot S(q^2_R)]=1\ot 1, \mbox{${\;\;\;}$}
\sum [1\ot S^{-1}(p^2_R)]q_R\Delta (p^1_R)=1\ot 1,
\end{equation}
\begin{equation}\label{pql}
\sum [S(p^1_L)\ot 1]q_L\Delta (p^2_L)=1\ot 1,%
\mbox{${\;\;\;}$}%
\sum \Delta (q^2_L)p_L[\smi (q^1_L)\ot 1]=1\ot 1,%
\end{equation}
\begin{eqnarray}
&&\hspace*{-2cm}(q_R\ot 1)(\Delta \ot id)(q_R)\Phi
^{-1}\nonumber\\
&=&\sum [1\ot S^{-1}(X^3)\ot S^{-1}(X^2)] [1\ot S^{-1}(f^2)\ot
S^{-1}(f^1)] (id \ot \Delta )(q_R\Delta
(X^1)),\label{qr2}\\
&&\hspace*{-2cm}\Phi (\Delta \ot
id)(p_R)(p_R\ot id)\nonumber\\
&=&\sum (id\ot \Delta )(\Delta (x^1)p_R)(1\ot f^{-1})(1\ot
S(x^3)\ot S(x^2))\label{pr1},
\end{eqnarray}
where $f=\sum f^1\ot f^2$ is the twist defined in (\ref{f}).

\subsection*{The smash product}
Suppose that $(H, \Delta , \varepsilon , \Phi )$ is a
quasi-bialgebra. If $U,V,W$ are left (right) $H$-modules, define
$a_{U,V,W}, {\bf a}_{U, V, W} :(U\otimes V)\otimes W\rightarrow
U\otimes (V\otimes W)$
by %
\begin{eqnarray*}
&&\hspace*{-2cm}a_{U,V,W}((u\otimes v)\otimes w)=\Phi \cdot (u\otimes
(v\otimes w)),\\[1mm]
&&\hspace*{-2cm}{\bf a}_{U, V, W}((u\ot v)\ot w)= (u\ot (v\ot w))\cd \Phi
^{-1}.
\end{eqnarray*}
Then the category $_H{\cal M}$ (${\cal M}_H$) of
left (right) $H$-modules becomes a monoidal category (see
\cite{k}, \cite{m2} for the terminology) with tensor product
$\otimes $ given via $\Delta $, associativity constraints
$a_{U,V,W}$ (${\bf a}_{U, V, W}$), unit $k$ as a trivial
$H$-module and the usual left and right
unit constraints.\\
Now, let $H$ be a quasi-bialgebra. We say that a $k$-vector space
$A$ is a left $H$-module algebra if it is an algebra in the
monoidal category $_H{\cal M}$, that is $A$ has a multiplication
and a usual unit $1_A$ satisfying the
following conditions:
\begin{eqnarray}
&&\hspace*{-2cm}(a a')a''=\sum (X^1\cd a)[(X^2\cd a')(X^3\cd
a'')],\label{ma1}\\[1mm]
&&\hspace*{-2cm}h\cd (a a')=\sum (h_1\cd a)(h_2\cd a'),
\label{ma2}\\[1mm]
&&\hspace*{-2cm}h\cd 1_A=\va (h)1_A\label{ma3}
\end{eqnarray}
for all $a, a', a''\in A$ and $h\in H$,
where $h\ot a\ra h\cd a$ is the $H$-module structure of $A$.
Following \cite{bpv} we define the smash product $A\# H$ as
follows: as a vector space $A\# H$ is $A\ot H$ ($a\ot h$
viewed as an element of $A\# H$
will be written $a\# h$) with multiplication
given by
\begin{equation}\label{sm1}
(a\# h)(a'\# h')=
\sum (x^1\cd a)(x^2h_1\cd a')\# x^3h_2h',
\end{equation}
for all $a, a'\in A$, $h, h'\in H$. $A\# H$ is an
associative algebra and it is defined by a universal property (as
Heyneman and Sweedler did for Hopf algebras, see \cite{bpv}). It
is easy to see that $H$ is a subalgebra of $A\# H$ via $h\mapsto
1\# h$, $A$ is a $k$-subspace of $A\# H$ via
$a\mapsto a\# 1$ and the following relations hold:
\begin{equation}\label{sm2}
(a\# h)(1\# h')=a\# hh', \mbox{${\;\;\;}$} (1\# h)(a\#
h')=\sum h_1\cd a\# h_2h', %
\end{equation}
for all $a\in A$, $h, h'\in H$.\\
For further use we need also the notion of right $H$-module
coalgebra. Suppose that $H$ is a quasi-bialgebra. Since the
category of right $H$-modules is a monoidal category we can define
coalgebras in this category. So, we say that a $k$-linear space
$C$ is a right $H$-module coalgebra if $C$ is a coalgebra in the
monoidal category ${\cal M}_H$, that is, if $C$ is a right
$H$-module (denote by $c\cd h$ the action of $h$ on $c$) and has a
comultiplication $\und :\ C\ra C\ot C$ and a usual counit $\une
:\ C\ra k$ satisfying the following relations
\begin{eqnarray}
&&(\und \ot id_C)(\und (c))\Phi ^{-1}=(id_C\ot \und )(\und (c))
\mbox{${\;\;}$$\forall $ $c\in C$},\label{rmc1}\\[1mm]%
&&\und (c\cd h)=\sum \una \cd h_1\ot \unb \cd h_2
\mbox{${\;\;}$$\forall $ $c\in C$, $h\in H$},
\label{rmc2}\\[1mm]
&&\une (c\cd h)=\une (c)\va (h) \mbox{${\;\;}$$\forall $ $c\in C$,
$h\in H$}\label{rmc3}
\end{eqnarray}
where we use the Sweedler-type notation
$$
\und (c)=\una \ot
\unb , \mbox{${\;\;}$} (\und \ot id_C)(\und (c))= \sum \unaa \ot
\unab \ot \unb \mbox{${\;\;}$ etc.}%
$$
%%%%%%%%%%%%%%%%%%%%%%%%%%%%%%%%%%%%%%%%%%%%%%%%%%%%%%
\section{The quasi-smash product}\selabel{2}
%%%%%%%%%%%%%%%%%%%%%%%%%%%%%%%%%%%%%%%%%%%%%%%%%%%%%%
\setcounter{equation}{0}
Recall from \cite{hn1} the notion of comodule algebra over a
quasi-bialgebra.
\begin{definition}
Let $H$ be a quasi-bialgebra. A unital associative algebra
$\mathfrak{A}$ is called a right $H$-comodule algebra if there
exist an algebra morphism $\r :\mathfrak{A}\ra \mathfrak{A}\ot H$
and an invertible element $\Phi_{\r }\in \mathfrak{A}\ot H\ot H$
such that
\begin{eqnarray}
&&\Phi _{\r }(\r \ot id)(\r (\mf {a}))=(id\ot \Delta
)(\r (\mf {a}))\Phi _{\r }
\mbox{${\;\;\;}$$\forall $ $\mf {a}\in
\mathfrak{A}$,}\label{rca1}\\[1mm]
&&(1_{\mf {A}}\ot \Phi)(id\ot \Delta \ot id)(\Phi _{\r })(\Phi
_{\r }\ot 1_H)= (id\ot id\ot \Delta )(\Phi _{\r })(\r \ot id\ot
id)(\Phi _{\r }),\label{rca2}\\[1mm]
&&(id\ot \va)\circ \r =id ,\label{rca3}\\[1mm]
&&(id\ot \va \ot id)(\Phi _{\r })=(id\ot id\ot \va )(\Phi _{\r }
)=1_{\mathfrak{A}}\ot 1_H.\label{rca4}
\end{eqnarray}
Similarly, a unital associative algebra $\mathfrak{B}$ is called
a left $H$-comodule algebra if there exist an algebra morphism $\l
: \mf {B}\ra H\ot \mathfrak{B}$ and an invertible element $\Phi
_{\l }\in H\ot H\ot \mathfrak{B}$ such that the following
relations hold
\begin{eqnarray}
&&(id\ot \l )(\l (\mf {b}))\Phi _{\l }=\Phi _{\l
}(\Delta \ot id)(\l (\mf {b}))%
\mbox{${\;\;\;}$$\forall $ $\mf {b}\in \mathfrak{B}$,}
\label{lca1}\\[1mm]%
&&(1_H\ot \Phi _{\l })(id\ot \Delta \ot id)(\Phi _{\l })(\Phi \ot
1_{\mf {B}})= (id\ot id\ot \l )(\Phi _{\l })(\Delta \ot id\ot
id)(\Phi _{\l }),\label{lca2}\\[1mm]%
&&(\va \ot id)\circ \l =id ,\label{lca3}\\[1mm]%
&&(id\ot \va \ot id)(\Phi _{\l })=(\va \ot id\ot id)(\Phi _{\l }
)=1_H\ot 1_{\mathfrak{B}}.\label{lca4}
\end{eqnarray}
\end{definition}

When $H$ is a quasi-bialgebra, particular examples of left and
right $H$-comodule algebras are given by $\mf {A}=\mf {B}=H$ and
$\r =\l =\Delta $,
$\Phi _{\r }=\Phi _{\l }=\Phi $.\\
For a right $H$-comodule algebra $({\mf A}, \r , \Phi _{\r })$ we
will denote
$$
\r (\mfa )=\sum \mfa _{<0>}\ot \mfa _{<1>}, \mbox{${\;\;}$} (\r
\ot id)(\r (\mfa ))=\sum \mfa _{<0, 0>}\ot \mfa _{<0, 1>} \ot \mfa
_{<1>} \mbox{${\;\;}$etc.}
$$
for any $\mfa \in {\mf A}$. Similarly, for a left $H$-comodule
algebra $({\mf B}, \l , \Phi _{\l })$, if $\mfb \in {\mf B}$ then
we will denote
$$
\l (\mfb )=\sum \mfb _{[-1]}\ot \mfb _{[1]}, \mbox{${\;\;}$}
(id\ot \l )(\l (\mfb ))=\sum \mfb _{[-1]}\ot \mfb _{[0,-1]}\ot
\mfb _{[0, 0]} \mbox{${\;\;}$etc.}
$$
In analogy with the notation of the reassociator $\Phi
$ of $H$, we will write
$$
\Phi _{\r }=\sum \tilde {X}^1_{\r }\ot \tilde {X}^2_{\r }\ot
\tilde {X}^3_{\r }=
\sum \tilde {Y}^1_{\r }\ot \tilde {Y}^2_{\r }\ot \tilde {Y}^3_{\r }={\rm~etc.}
$$
and its inverse by
$$
\Phi _{\r }^{-1}=\sum \tilde {x}^1_{\r }\ot \tilde {x}^2_{\r }\ot
\tilde {x}^3_{\r }=\sum \tilde {y}^1_{\r }\ot \tilde {y}^2_{\r
}\ot \tilde {y}^3_{\r }={\rm~etc.}
$$
Similarly for the element $\Phi _{\l }$ of a left $H$-comodule
algebra $\mf {B}$.  If there is no danger of confusion we will omit
the subscription $\r $ or $\l $ for the tensor components of the
elements $\Phi _{\r }$, $\Phi _{\l }$ or for the tensor components
of the elements  $\Phi _{\r }^{-1}$, $\Phi
_{\l }^{-1}$.\\
Suppose now that $H$ is finite dimensional and $\mf {A}$ is a
right $H$-comodule algebra. In the Hopf case, we have that $\mf
{A}$ is a left $H^*$-module algebra so we can consider the smash
product $\mf {A}\# H^*$. We will see that a similar result holds
in the quasi-Hopf case. Since the resulting object of $\mf {A}$ and
$H^*$ is an algebra in the monoidal category $_H{\cal M}$ (i.e. an
$H$-module algebra) we will call it the quasi-smash product
between $\mf {A}$ and $H^*$. As expected, when $H$ is a bialgebra
the two constructions, the smash product and the quasi-smash
product, coincide.\\
Let $\{e_i\}_{i=\ov {1, n}}$ be a basis of $H$, and
$\{e^i\}_{i=\ov {1, n}}$ the corresponding dual basis of $H^*$.
$H^*$ is a coassociative coalgebra, with comultiplication %
$$
\widehat {\Delta }(\v )= \sum \v _1\ot \v _2= \sum_{i,
j=1}^n\v (e_ie_j)e^i\ot e^j,
$$
or, equivalently,
$$
\widehat {\Delta }(\v )=\sum \v _1\ot \v_2~~ \Longleftrightarrow ~~\v
(hh')=\sum \v _1(h)\v _2(h'),%
\mbox{${\;\;\;}$$\forall h, h'\in H$.}
$$
$H^*$ is also an $(H,H)$-bimodule, by
$$
<h\rh \v , h'>=\v (h'h), \mbox{${\;\;\;}$}
<\v \lh h, h'>=\v (hh').
$$
The convolution $<\v \psi , h>=\sum \v (h_1)\psi (h_2)$, $h\in H$,
is a multiplication on $H^*$; it is not associative,
but only quasi-associative:
\begin{equation}\label{mbia1}
[\v \psi]\xi=\sum (X^1\rh \v \lh x^1)[(X^2\rh \psi \lh x^2)
(X^3\rh \xi \lh x^3)], \mbox{${\;\;\;}$$\forall \v , \psi , \xi
\in H^*$.}
\end{equation}
In addition, for all $h\in H$ and $\v , \psi \in H^*$ we have
that
\begin{equation}\label{mbia2}
h\rh (\v \psi )=\sum (h_1\rh \v )(h_2\rh \psi )
\mbox{${\;\;\;}$and${\;\;\;}$}
(\v \psi )\lh h=\sum (\v \lh h_1)(\psi \lh h_2).
\end{equation}
In other words, $H^*$ is an algebra in the monoidal category of
$(H,H)$-bimodules $_H{\cal M}_H$.\\
Now, if $(\mf {A}, \r , \Psi _{\r })$ is a right $H$-comodule
algebra we define a multiplication on $\mf {A}\ot H^*$ as follows
\begin{equation}\label{qsm}
(\mf {a}\ovsm \v )(\mf {a}'\ovsm \Psi )= \sum \mf {a}\mf
{a}'_{<0>}\tx ^1\ovsm (\v \lh \mf {a}'_{<1>}\tx ^2)
(\psi \lh \tx ^3)
\mbox{${\;\;\;}$$\forall $$\mf {a}, \mf {a}'\in \mf {A}$ and
$\v , \Psi \in H^*$}
\end{equation}
where we write $\mfa \ovsm \v $ for $\mfa \ot \v $, $\r (\mfa
)=\sum \mfa _{<0>}\ot \mfa _{<1>}$, and $\Phi _{\r }^{-1}=\sum
\tx ^1\ot \tx ^2\ot \tx ^3$. We denote this structure on $\mf
{A}\ot H^*$ by $\mf {A}\ovsm H^*$.

\begin{proposition}
Let $H$ be a finite dimensional quasi-bialgebra and
$(\mf {A}, \r, \Phi _{\r })$ be a right $H$-comodule algebra. Then $\mf
{A}\ovsm H^*$ is a $H$-module algebra with unit $1_{\mf {A}}\ovsm
\va $ and with the left $H$-action given by
\begin{equation}\label{aqsm}
h\cd (\mf {a}\ovsm \v )=\mf {a}\ovsm h\rh \v,~~~
\forall ~h\in H,~\mf {a}\in \mf {A}, \hbox{ and}~\v
\in H^*.
\end{equation}
\end{proposition}

\begin{proof}
Since $H^*$ is a left $H$-module via the action $\rh $, it is easy
to see that $\mf {A}\ovsm H^*$ is a left $H$-module via the action
(\ref{aqsm}). Now, we will prove that $\mf {A}\ovsm
H^*$ is an algebra in $_H{\cal M}$ with unit $1_{\mf {A}}\ovsm \va
$. Indeed, for all $\mfa , \mfa', \mfa''\in \mf {A}$ and
$\v , \psi , \chi \in H^*$
\begin{eqnarray*}
&&\hspace*{-2cm}
[X^1\cd (\mfa \ovsm \v )]\{[X^2\cd (\mfa'\ovsm \psi )][X^3\cd
(\mfa''\ovsm \chi )]\}\\
&=&\sum (\mfa \ovsm X^1\rh \v )[(\mfa'\ovsm X^2\rh \psi )(\mfa''\ovsm X^3\rh
\chi )]\\
&=&\sum (\mfa \ovsm X^1\rh \v )[\mfa'\mfa''_{<0>}\tx
^1\ovsm (X^2\rh \psi \lh \mfa''_{<1>}\tx ^2)(X^3\rh \chi \lh
\tx ^3)]\\
{\rm (\ref{mbia2})}&=&\sum \mfa \mfa'_{<0>}\mfa''_{<0,
0>}\tx ^1_{<0>}\ty ^1\ovsm (X^1\rh \v \lh \mfa'_{<1>}\mfa''_{<0, 1>}\tx
^1_{<1>}\ty ^2)\\
&&[(X^2\rh \psi \lh \mfa''_{<1>}\tx ^2\ty ^3_1)(X^3\rh \chi \lh \tx ^3\ty
^3_2)]\\
{\rm (\ref{mbia1},\ref{rca2})}&=&\sum \mfa \mfa'_{<0>}\mfa''_{<0, 0>}\tx
^1\ty ^1\ovsm [(\v \lh \mfa'_{<1>}\mfa''_{<0, 1>}\tx ^2\ty ^2_1)(\psi \lh
\mfa''_{<1>}\tx ^3\ty ^2_2)]\\
&&(\chi \lh \ty ^3)\\
{\rm (\ref{rca1},\ref{mbia2})}&=&\sum \mfa \mfa
'_{<0>}\tx ^1\mfa''_{<0>}\ty ^1\ovsm \{[(\v \lh \mfa
'_{<1>}\tx ^2)(\psi \lh
\tx ^3)]\lh \mfa''_{<1>}\ty ^2\}(\chi \lh \ty ^3)\\
&=&\sum [\mfa \mfa'_{<0>}\tx ^1\ovsm (\v \lh \mfa'_{<1>}\tx ^2)(\psi \lh \tx
^3)](\mfa''\ovsm \chi )\\
&=&[(\mfa \ovsm \v )(\mfa'\ovsm \psi )](\mfa''\ovsm \chi ).
\end{eqnarray*}
It is not hard to see that $1_{\mf {A}}\ovsm \va $ is the unit of
$\mf {A}\ovsm H^*$ and that $h\cd (1_{\mf {A}}\ovsm \va)=\va
(h)1_{\mf {A}}\ovsm \va $ for all $h\in H$. Finally, for all $h\in
H$, $\mfa , \mfa'\in \mf {A}$ and $\v , \psi \in H^*$, we
calculate:
\begin{eqnarray*}
&&\hspace*{-2cm}
\sum [h_1\cd (\mfa \ovsm \v )][h_2\cd (\mfa'\ovsm \psi )]
=\sum (\mfa \ovsm h_1\rh \v )(\mfa'\ovsm h_2\rh \psi )\\
&=&\sum \mfa \mfa'_{<0>}\tx ^1\ovsm (h_1\rh \v \lh \mfa
'_{<1>}\tx ^2)(h_2\rh \psi \lh \tx ^3)\\
{\rm (\ref{mbia2})}&=&\sum \mfa \mfa '_{<0>}\tx ^1\ovsm
h\rh [(\v \lh \mfa'_{<1>}\tx ^2)(\psi \lh \tx ^3)]\\
{\rm (\ref{aqsm})}&=&h\cd [(\mfa \ovsm \v )(\mfa'\ovsm
\psi )].
\end{eqnarray*}
\end{proof}

$(H, \Delta , \Phi )$ is a right $H$-comodule algebra, so it makes
sense to consider the quasi-smash product $H\ovsm H^*$. In this case
where $H$ is a Hopf algebra, $H\# H^*$ is called the Heisenberg double of
$H$, and we will keep the same terminology for quasi-Hopf algebras.
${\cal H}(H)=H\ovsm H^*$ is
not an associative algebra but it is an algebra in the monoidal
category $_H{\cal M}$. If $H$ is a finite dimensional Hopf algebra
then ${\cal H}(H)$ is isomorphic to the algebra $\End_k(H)$.
In order to prove a similar result for a finite dimensional
quasi-Hopf algebra, we first have
to deform the algebra structure of $\End_k(H)$.

\begin{proposition}\prlabel{2.3}
Let $H$ be a finite dimensional quasi-Hopf algebra. Define
$$
\mu :\ H\ovsm H^*\ra \End_k(H),
\mbox{${\;\;\;}$}
\mu (h\ovsm \v )(h')=\sum \v (h'_2p^2_L)hh'_1p^1_L$$
for all $h, h'\in H$ and $\v \in H^*$,
where $p_L=\sum p^1_L\ot p^2_L$ is the element defined by
(\ref{ql}). Then $\mu $ is a bijection, and therefore there exists a
unique $H$-module algebra structure on $\End_k(H)$ such that $\mu $
becomes an $H$-module algebra isomorphism. The multiplication, the
unit and the $H$-module structure of $\End_k(H)$ are given by
\begin{eqnarray}
&&(u\ov {\circ }v)(h)=\sum u(v(hx^3X^3_2)\smi (S(x^1X^2)\a
x^2X^3_1))\smi (X^1)\label{end}\\
&&1_{\End_k(H)}=\smi (\b )\rh id_H~~~;~~~(h\cdot u)(h')=u(h'h_2)\smi (h_1)
\label{end1}
\end{eqnarray}
for all $u, v\in \End_k(H)$ and $h, h'\in H$.
\end{proposition}

\begin{proof}
Let $\{e_i\}_{i=\ov {1, n}}$ be a basis of $H$ and $\{e^i\}_{i=\ov
{1, n}}$ the corresponding dual basis of $H^*$. We claim that the
inverse of $\mu $ is $\mu ^{-1}:\ \End_k(H)\ra
H\ovsm H^*$ given by
$$
\mu ^{-1}(u)=\sum u(q^2_L(e_i)_2)\smi (q^1_L(e_i)_1)\ovsm e^i$$
for all $u\in \End_k(H)$,
where $q_L=\sum q^1_L\ot q^2_L$ is the element defined by
(\ref{ql}). Indeed,
for any $h\in H$ and $\v \in H^*$ we have:
\begin{eqnarray*}
&&\hspace*{-2cm}
(\mu ^{-1}\circ \mu )(h\ovsm \v )
=\sum_{i=1}^n\mu (h\ovsm \v )(q^2_L(e_i)_2)\smi
(q^1_L(e_i)_1)\ovsm e^i\\
&=&\sum_{i=1}^n\v ((q^2_L)_2(e_i)_{(2,
2)}p^2_L)h(q^2_L)_1(e_i)_{(2, 1)}p^1_L\smi (q^1_L(e_i)_1)\ovsm e^i\\
{\rm (\ref{ql1})}&=&\sum_{i=1}^n\v
((q^2_L)_2p^2_Le_i)h(q^2_L)_1p^1_L\smi (q^1_L)\ovsm e^i\\
{\rm (\ref{pql})}&=&\sum_{i=1}^n\v (e_i)h\ovsm
e^i=h\ovsm \v
\end{eqnarray*}
and, in a similar way, for $u \in \End_k(H)$ and $h\in H$ we have that
\begin{eqnarray*}
&&\hspace*{-2cm}
(\mu \circ \mu ^{-1})(u)(h)
=\sum_{i=1}^n\mu (u(q^2_L(e_i)_2)\smi
(q^1_L(e_i)_1)\ovsm e^i)(h)\\
&=&\sum_{i=1}^ne^i(h_2p^2_L)u(q^2_L(e_i)_2)\smi
(q^1_L(e_i)_1)h_1p^1_L\\
&=&\sum u(q^2_Lh_{(2, 2)}(p^2_L)_2)\smi (q^1_Lh_{(2,
1)}(p^2_L)_1)h_1p^1_L\\
{\rm (\ref{ql1})}&=&\sum u(hq^2_L(p^2_L)_2)\smi
(q^1_L(p^2_L)_1)p^1_L\\
{\rm (\ref{pql})}&=&u (h)1_H=u(h).
\end{eqnarray*}
Using the bijection $\mu$, we transport the $H$-module algebra
structure from $H\ovsm H^*$ to $\End_k(H)$. First we compute the transported
multiplication $\ov {\circ }$: for all $u, v\in
\End_k(H)$, we find
\begin{eqnarray*}
u\ov {\circ }v
&=&\sum \mu (\mu ^{-1}(u)\mu ^{-1}(v))\\
&=&\sum_{i, j=1}^n\mu ((u(q^2_L(e_i)_2)\smi
(q^1_L(e_i)_1)\ovsm
e^i)(v(Q^2_L(e_j)_2)\smi (Q^1_L(e_j)_1)\ovsm e^j))\\
{\rm (\ref{qsm})}&=&\sum_{i, j=1}^n
\mu (u(q^2_L(e_i)_2)\smi (q^1_L(e_i)_1)
[v(Q^2_L(e_j)_2)\smi (Q^1_L(e_j)_1)]_1x^1\\
&&\ovsm (e^i\lh [v(Q^2_L(e_j)_2)\smi (Q^1_L(e_j)_1)]_2x^2)(e^j\lh
x^3))
\end{eqnarray*}
where $\sum Q^1_L\ot Q^2_L$ is another copy of $q_L$.
Note that (\ref{q3}) and (\ref{ql}) imply that
\begin{equation}\label{lq}
\sum S(x^1)q^1_Lx^2_1\ot q^2_Lx^2_2\ot x^3=\sum q^1_LX^1\ot
(q^2_L)_1X^2\ot (q^2_L)_2X^3.
\end{equation}
Hence, for all $h\in H$ we have that
\begin{eqnarray*}
(u\ov {\circ }v)(h)&\pile{{\rm (\ref{q1})}\\ =}&
\sum_{i, j, l=1}^n<e^l,
v(Q^2_L(e_j)_2)\smi (Q^1_L(e_j)_1)h_1><e^i,
(e_l)_2x^2(p^2_L)_1>\\%
&&<e^j, h_2x^3(p^2_L)_2>
u(q^2_L(e_i)_2)\smi (q^1_L(e_i)_1)(e_l)_1x^1p^1_L\\%
&=&\sum_{i, l=1}^n<e^l, v(Q^2_Lh_{(2,
2)}x^3_2(p^2_L)_{(2, 2)})\smi (Q^1_Lh_{(2, 1)}x^3_1(p^2_L)_{(2,
1)})h_1>\\%
&&<e^i, (e_l)_2x^2(p^2_L)_1>
u(q^2_L(e_i)_2)\smi (q^1_L(e_i)_1)(e_l)_1x^1p^1_L\\%
{\rm (\ref{ql1})}&=&\sum_{l=1}^n <e^l,
v(hQ^2_Lx^3_2(p^2_L)_{(2, 2)})\smi (Q^1_Lx^3_1(p^2_L)_{(2, 1)})>
u(q^2_L(e_l)_{(2, 2)}(x^2(p^2_L)_1)_2)\\%
&&\smi (q^1_L(e_l)_{(2, 1)}(x^2(p^2_L)_1)_1)(e_l)_1x^1p^1_L\\%
{\rm (\ref{ql1}), (\ref{lq})}&=&\sum_{l=1}^n <e^l,
v(hQ^2_L(q^2_L)_{(2,2)}X^3_2((p^2_L)_2)_2)\smi (Q^1_L(q^2_L)_{(2,
1)}X^3_1((p^2_L)_2)_1)>\\%
&&u(e_l(q^2_L)_1X^2(p^2_L)_{(1, 2)})\smi (q^1_LX^1(p^2_L)_{(1,
1)})p^1_L\\%
{\rm (\ref{q1})}&=&\sum u(v(hQ^2_L[q^2_L(p^2_L)_2]_{(2,
2)}X^3_2)\smi (Q^1_L([q^2_L(p^2_L)_2]_{(2,
1)}X^3_1)[q^2_L(p^2_L)_2]_1X^2)\\%
&&\smi (q^1_L(p^2_L)_1X^1)p^1_L\\%
{\rm (\ref{pql})}&=&\sum u(v(hQ^2_LX^3_2)\smi
(Q^1_LX^3_1)X^2)\smi (X^1)\\%
&=&\sum u(v(hx^3X^3_2)\smi (S(x^1X^2)\a x^2X^3_1))\smi (X^1).
\end{eqnarray*}
Thus, we have obtained (\ref{end}). Similar
computations show that the tranported unit and the $H$-action on
$\End_k(H)$ are given by (\ref{end1}).
\end{proof}

\begin{remarks}\relabel{2.4}
Let $H$ be a finite dimensional quasi-Hopf algebra,
$\{e_i\}_{i=\ov {1, n}}$ a basis of $H$ and $\{e^i\}_{i=\ov
{1, n}}$ the corresponding dual basis of $H^*$.\\
1) The bijection $\mu $ defined in \prref{2.3} induces an
associative algebra structure on the $k$-vector space $H\ot H^*$:
it suffices to transport the composition on $\End_k(H)$ to
$H\ot H^*$.
If $H$ is a Hopf algebra, then we find the smash product, and
the corresponding category of left or right representations
is equivalent to a certain category of Hopf modules ${\cal M}_H^H$.
This could open the door to defining the category of Hopf modules
over a quasi-Hopf algebra. Unfortunately, the involved structures
are not easy to describe, except in the situations where
$H$ is a twisted Hopf algebra. A possible way to define ${\cal M}_H^H$ will be
presented in the next Section.\\
2) Let $(\mf {A}, \r , \Phi _{\r })$ be a right $H$-comodule
algebra. As in the Hopf case, it is possible to associate
different (quasi)smash products to ${\mf A}$. Observe first that
the map $\nu :\ {\mf A}\ovsm H^*\ra \Hom_k(H, {\mf A})$ given by
$\nu (\mfa \ovsm \v )(h)=\v (h)\mfa $, for all $\mfa \in {\mf A}$,
$\v \in H^*$ and $h\in H$, is a $k$-linear isomorphism. The
inverse of $\nu $ is given by the formula%
$$%
\nu ^{-1}(w)=\sum \limits_{i=1}^nw(e_i)\ovsm e^i%
$$%
for $w\in \Hom_k(H, {\mf A})$. Secondly, by
transporting the quasi-smash algebra structure from ${\mf A}\ovsm
H^*$ to $\Hom_k(H, {\mf A})$ via the isomorphism $\nu $, we obtain
that $\Hom_k(H, {\mf A})$ is an $H$-module algebra. So, if $H$ is an
arbitrary quasi-Hopf algebra and $(\mf {A}, \r , \Phi _{\r })$ is
a right $H$-comodule algebra, then we can define the quasi-smash
product $\ovsm (H, {\mf A})$ as follows: $\ovsm (H, {\mf A})$ is
the $k$-vector space $Hom_k(H, {\mf A})$ with multiplication given
by
\begin{equation}\label{sovsm}
(v*w)(h)=\sum v(w(\tx ^3h_2)_{<1>}\tx ^2h_1)w(\tx ^3h_2)_{<0>}\tx
^1
\end{equation}
for $v, w\in \ovsm (H, {\mf A})$ and $h\in H$.
The unit is $1_{\ovsm (H, {\mf A})}(h)=\va (h)1_{\mf A}$ and the
$H$-module structure is given by $(h\cd v)(h')=v(h'h)$, $h,
h'\in H$, $v\in \Hom_k(H, {\mf A})$. Of course, if $H$ is finite
dimensional then ${\mf A}\ovsm H^*\simeq \ovsm (H, {\mf A})$ as
$H$-module algebras.
\end{remarks}

%%%%%%%%%%%%%%%%%%%%%%%%%%%%%%%%%%%%%%%%%%%%%%%%%%%%%%%
\section{Two-sided Hopf modules and relative Hopf
modules}\selabel{3}
%%%%%%%%%%%%%%%%%%%%%%%%%%%%%%%%%%%%%%%%%%%%%%%%%%%%%%%%%%
\setcounter{equation}{0}%
We cannot define comodules over a
quasi-bialgebra $H$, because a quasi-bialgebra is not
coassociative. Also the definition of Hopf modules is quite
complicated. A possible approach is the following: assume that we
can deform the comultiplication on $H$ in such a way that $H$ with
this new comultiplication (denoted $\ov{H}$) becomes a coalgebra
in the monoidal category ${\cal M}_H$, and define ${\cal M}_H^H$
as the category of right $[\ov {H}, H]$-Hopf modules, as
introduced in \cite{bn}. More precisely, a right $H$-Hopf module
is a right $H$-module which is also a right $\ov {H}$-comodule in
the monoidal category ${\cal M}_H$ (for more detail, see
\seref{5}). This machinere works if $H$ is a twisted bialgebra.
Indeed, if $H$ is an ordinary bialgebra, and $F\in H\ot H$ is a
twist, then $H_F$ is a quasi-bialgebra with counit $\va $ and
comultiplication given by (\ref{g1}). If we define on $H_F$ a new
comultiplication given by
$$%
\ov {\Delta }_F(h)=\sum h_1{\cal G}^1\ot h_2{\cal G}^2%
$$%
for all $h\in H$, where $F^{-1}=\sum {\cal G}^1\ot {\cal G}^2$,
and we denote $H_F$ with this new comultiplication by $\ov {H}_F$,
then it is not hard to see that $\ov {H}_F$ is a coalgebra in the
monoidal category ${\cal M}_{H_F}$. Therefore, by a right
$H_F$-Hopf module, we mean a $k$-vector space $M$, which is a
right $H$-module (the right $H$-action of $h\in H$ on $m\in M$ is
denoted by $mh$), together with a $k$-linear map $\r_M:\ M\ra M\ot
H$ such that the following relations hold, for all $m\in M$ and
$h\in H$:
\begin{eqnarray*}
&&(\r _M\ot id_H)(\r _M(m))(F\ot 1_H)=(id_M\ot \Delta )(\r
_M(m)F),\\%
&&(id_M\ot \va )(\r _M(m))=m,\\%
&&\r _M(mh)F=\r _M(m)F\Delta (h).
\end{eqnarray*}
We can conclude that categories of Hopf modules over quasi-Hopf
algebras can be defined using  (co)algebras and (co)modules in
monoidal categories. This point of view was used in \cite{bn},
\cite{hn3} and \cite{sh} in order to define the categories of
relative Hopf modules, quasi-Hopf bimodules and
two-sided two-cosided Hopf modules. In the sequel, we will study
all these categories in a more general
context.

\begin{definition}\delabel{3.1}
Let $H$ be a quasi-Hopf algebra and $({\mf A}, \r , \Phi _{\r })$
a right $H$-comodule algebra. A two-sided $({\mf A}, H)$-Hopf module
is an $({\mf A}, H)$-bimodule $M$ together with a $k$-linear map
$$\r _M:\ M\ra M\ot H,~~\r _M(m)=\sum m_{(0)}\ot m_{(1)}$$
satisfying the following relations, for all $m\in M$, $h\in H$ and $\mfa
\in {\mf A}$.
The actions of $h\in H$ and $\mfa \in {\mf A}$ on $m\in M$ are
denoted by $h\succ m$ amd $m\prec \mfa$.
\begin{eqnarray}
&&\hspace{-1cm}(id_M\ot \va )\circ \r _M=id_M,\label{thm1}\\
&&\hspace{-1cm}\Phi \cdot (\r _M\ot id_H)(\r _M(m))=(id_M\ot
\Delta )(\r _M
(m))\cdot \Phi _{\r }, \label{thm2}\\%
&&\hspace{-1cm}\r _M(h\succ m)=\sum h_1\succ m_{(0)}\ot
h_2m_{(1)},\label{thm3}\\%
&&\hspace{-1cm} \r _M(m\prec \mfa )=\sum m_{(0)}\prec \mfa
_{<0>}\ot m_{(1)}\mfa _{<1>}.\label{thm3a}
\end{eqnarray}
\end{definition}

The category of two-sided $({\mf A}, H)$-Hopf modules and
left $H$-linear, right ${\mf A}$-linear and right $H$-colinear maps
is denoted by ${}_H{\cal M}^H_{\mf A}$.\\
Observe that the category of two-sided $(H,
H)$-Hopf bimodules is nothing else then the category of right quasi-Hopf
$H$-bimodules introduced in \cite{hn3}.\\
We will use the following notation, similar to the notation for the
comultiplication on a quasi-bialgebra:
\begin{eqnarray*}
(\r _M\ot id_H)(\r _M(m))&=&\sum m_{(0, 0)}\ot m_{(0, 1)}\ot m_{(1)},\\%
(id_M\ot \Delta_H)(\r _M(m))&=&\sum m_{(0)}\ot m_{(1)_1}\ot
m_{(1)_2}.
\end{eqnarray*}

\begin{examples}\exlabel{3.2}
Let $H$ be a quasi-Hopf algebra and $({\mf A}, \r , \Phi_{\r })$
a right $H$-comodule algebra.\\

1) ${\cal V}={\mf A}\ot H\in {}_H{\cal M}^H_{\mf A}$. The structure maps
are
$$
h\succ (\mfa \ot h')=\mfa \ot hh'~~~;~~~
(\mfa \ot h)\prec \mfa'=\sum \mfa \mfa'_{<0>}\ot h\mfa'_{<1>}
$$
and
$$
\r_{\cal V}(\mfa \ot h)=\sum \mfa \tX ^1\ot h_1\tX ^2\ot h_2\tX ^3
$$
for all $h, h'\in H$ and $\mfa , \mfa '\in {\mf A}$. We leave
verification of the detail to the reader.\\%

2) ${\cal U}=H\ot {\mf A}\in {}_H{\cal M}^H_{\mf A}$. Now the structure
maps are given by the following formulas, for all
$h, h'\in H$ and $\mfa , \mfa '\in {\mf A}$:
$$h\succ (h'\ot \mfa )=hh'\ot \mfa~~~;~~~
(h\ot \mfa )\prec \mfa '=h\ot \mfa \mfa'$$
and
\begin{equation}\eqlabel{coaction}
\r _{\cal U}(h\ot \mfa )=\sum h_1\smi (q^2_L\tX ^3_2g^2)\ot \tX
^1\mfa_{<0>}\ot h_2\smi (q^1_L\tX ^3_1g^1)\tX ^2\mfa _{<1>}.
\end{equation}
Here $q_L=\sum q^1_L\ot q^2_L$ and $f^{-1}=\sum g^1\ot g^2$ are
the elements defined by the formulas (\ref{ql}) and (\ref{g}).\\%

Now consider $\theta :\ {\cal V}\ra {\cal U}$ given by
$$%
\theta (\mfa \ot h)=\sum h\smi (\mfa _{<1>}\tilde {p}^2_{\r })\ot
\mfa _{<0>}\tilde {p}^1_{\r }
$$%
for all $h\in H$ and $\mfa \in {\mf A}$, where we use the notation
\begin{equation}\label{tpr}
\tilde {p}_{\r }=\sum \tilde {p}^1_{\r }\ot \tilde {p}^2_{\r}
=\sum \tx ^1\ot \tx ^2\b S(\tx ^3)\in {\mf A}\ot H.
\end{equation}
$\theta $ is bijective; its inverse $\theta
^{-1}:\ {\cal U}\ra {\cal V}$ is defined as follows
$$
\theta ^{-1}(h\ot \mfa )=\sum \tilde {q}^1_{\r }\mfa _{<0>}\ot
h\tilde {q}^2_{\r }\mfa _{<1>}$$
with the notation
\begin{equation}\label{tqr}
\tilde {q}_{\r }=\sum \tilde {q}^1_{\r }\ot \tilde {q}^2_{\r}
=\sum \tX ^1\ot \smi (\a \tX ^3)\tX ^2\in {\mf A}\ot H.
\end{equation}
Furthermore, $\theta$ is a morphism of two-sided $({\mf A}, H)$-Hopf bimodules,
and we conclude that ${\cal U}=H\ot {\mf A}$ and ${\mf A}\ot H={\cal V}$
are isomorphic in $_H{\cal M}^H_{\mf A}$.\\
To prove this, we proceed as follows. First, by
\cite[Lemma 9.1]{hn1}, we have the following relations, for all $\mfa \in
{\mf A}$:
\begin{eqnarray}
&&\hspace*{-1cm}\sum \r (\mfa_{<0>})\tilde {p}_{\r }[1_{\mf A}\ot S(\mfa
_{<1>})]
=\tilde {p}_{\r }[\mfa \ot 1_H],\label{tpqr1}\\%
&&\hspace*{-1cm}\sum [1_{\mf A}\ot \smi (\mfa_{<1>})]\tilde {q}_{\r }\r
(\mfa_{<0>})=[\mfa \ot 1_H]\tilde {q}_{\r },\label{tpqr1a}\\%
&&\hspace{-1cm}\sum \r (\tilde {q}^1_{\r })\tilde {p}_{\r }[1_{\mf
A}\ot S(\tilde
{q}^2_{\r })]=1_{\mf A}\ot 1_H\label{tpqr2},\\%
&&\hspace*{-1cm}
\sum [1_{\mf A}\ot \smi (\tilde {p}^2_{\r })]\tilde {q}_{\r }\r
(\tilde {p}^1_{\r })=1_{\mf A}\ot 1_H,\label{tpqr2a}\\%
&&\hspace{-1cm}\Phi _{\r }(\r \ot id_H)(\tilde {p}_{\r })\tilde
{p}_{\r }= \sum (id \ot \Delta )(\r (\tx ^1)\tilde {p}_{\r
})(1_{\mf A}\ot
g^1S(\tx ^3)\ot g^2S(\tx ^2)),\label{tpr2}\\%
&&\hspace{-1cm}(\tilde {q}_{\r }\ot 1_H)(\r \ot id_H)(\tilde
{q}_{\r })\Phi _{\r }^{-1}=\sum [1_{\mf A}\ot \smi (f^2\tX ^3)\ot
\smi (f^1\tX ^2)](id _{\mf A}\ot \Delta )(\tilde {q}_{\r }\r (\tX
^1)).\label{tqr2}
\end{eqnarray}
Here $f=\sum f^1\ot f^2$ is the element defined in (\ref{f}) and
$f^{-1}=\sum g^1\ot g^2$. Using (\ref{tpqr1}-\ref{tpqr2a}) we can
show easily that $\theta $ and $\theta ^{-1}$ are inverses, and
that ${\cal U}$ is an $(H,{\mf A})$-bimodule via the actions
$\succ $ and $\prec $. We will finally compute the right
$H$-coaction on ${\cal U}$ transported from the coaction on ${\cal
V}$ using $\theta$, and then see that it coincides with
\eqref{coaction}. First observe that (\ref{tpr},\ref{rca2}) and
(\ref{rca4}) imply
\begin{equation}\label{tprr}
\sum \tX ^1_{<1>}\tilde {p}^2_{\r }S(\tX ^2)\ot \tX ^1_{<0>}\tilde
{p}^1_{\r }\ot \tX ^3=\sum \tx ^2S(\tx ^3_1p^1_L)\ot \tx ^1\ot \tx
^3_2p^2_L,
\end{equation}
where $p_L=\sum p^1_L\ot p^2_L$ is the element defined in
(\ref{ql}). Therefore, for all $h\in H$ and $\mfa \in {\mf A}$ we
have:
\begin{eqnarray*}
&&\hspace*{-2cm}(\theta \ot id_H)\circ \r _{\cal V}
\circ \theta ^{-1}(h\ot \mfa )\\
&=&\sum (\theta \ot id_H)(\tilde {q}^1_{\r }\mfa _{<0>}\tX ^1 \ot
h_1(\tilde {q}^2_{\r })_1\mfa _{<1>_1}\tX ^2\ot h_2(\tilde
{q}^2_{\r })_2\mfa _{<1>_2}\tX ^3)\\
&=&\sum h_1(\tilde {q}^2_{\r })_1\mfa _{<1>_1}\tX ^2\smi ((\tilde
{q}^1_{\r })_{<1>}\mfa _{<0, 1>}\tX ^1_{<1>}\tilde {p}^2_{\r })\ot
(\tilde {q}^1_{\r })_{<0>}\mfa _{<0, 0>}\tX
^1_{<0>}\tilde {p}^1_{\r }\\
&&\ot h_2(\tilde {q}^2_{\r })_2\mfa _{<1>_2}\tX ^3\\
{\rm (\ref{tprr})}&=&\sum h_1(\tilde {q}^2_{\r })_1\mfa
_{<1>_1}\tx ^3_1p^1_L\smi ((\tilde {q}^1_{\r })_{<1>}\mfa _{<0,
1>}\tx ^2)\ot (\tilde {q}^1_{\r })_{<0>}\mfa
_{<0, 0>}\tx ^1\\
&&\ot h_2(\tilde {q}^2_{\r })_2\mfa _{<1>_2}\tx ^3_2p^2_L\\
{\rm (\ref{rca1})}&=&\sum h_1(\tqrb )_1\tx ^3_1\mfa _{<1>_{(2,
1)}}p^1_L\smi ((\tqra )_{<1>}\tx ^2\mfa _{<1>_1})\ot (\tqra
)_{<0>}\tx ^1\mfa _{<0>}\\
&&\ot h_2(\tqrb )_2\tx ^3_2\mfa _{<1>_{(2, 2)}}p^2_L\\
{\rm (\ref{ql1})}&=&\sum h_1(\tqrb )_1\tx ^3_1p^1_L\smi
((\tqra )_{<1>}\tx ^2)\ot (\tqra )_{<0>}\tx ^1\mfa _{<0>}\ot
h_2(\tqrb )_2\tx ^3_2p^2_L\mfa _{<1>}\\
{\rm (\ref{tqr2})}&=&\sum h_1\smi (\ty ^3g^2)\tqrb (\tQra
)_{<1>}\tz ^2(\ty ^1_{<1>}\tx ^3)_1p^1_L\smi ((\tqra
(\tQra )_{<0>}\tz ^1\ty ^1_{<0>})_{<1>}\tx ^2)\\
&&\ot (\tqra (\tQra )_{<0>}\tz ^1\ty ^1_{<0>})_{<0>}\tx ^1\mfa
_{<0>}\ot h_2\smi (\ty ^2g^1)\tQrb \tz ^3(\ty ^1_{<1>}\tx
^3)_2p^2_L\mfa _{<1>}\\
{\rm (\ref{rca1})}&=&\sum h_1\smi (\ty ^3g^2)\tqrb (\tQra
)_{<1>}\tz ^2\tx ^3_1\ty ^1_{<1>_{(2, 1)}}p^1_L\smi ((\tqra (\tQra
)_{<0>})_{<1>}\tz ^1_{<1>}\tx ^2\ty ^1_{<1>_1})\\
&&\ot (\tqra (\tQra )_{<0>})_{<0>}\tz ^1_{<0>}\tx ^1\ty
^1_{<0>}\mfa _{<0>}\ot h_2\smi (\ty ^2g^1)\tQrb \tz ^3\tx ^3_2\ty
^1_{<1>_{(2, 2)}}p^2_L\mfa _{<1>}\\
{\rm (\ref{ql1})}&=&\sum h_1\smi (\ty ^3g^2)\tqrb (\tQra
)_{<1>}\tz ^2\tx ^3_1p^1_L\smi ((\tqra (\tQra
)_{<0>})_{<1>}\tz ^1_{<1>}\tx ^2)\\
&&\ot (\tqra (\tQra )_{<0>})_{<0>}\tz ^1_{<0>}\tx ^1\ty
^1_{<0>}\mfa _{<0>}\ot h_2\smi (\ty ^2g^1)\tQrb \tz ^3\tx
^3_2p^2_L\ty ^1_{<1>}\mfa _{<1>}\\
{\rm (\ref{ql},\ref{rca2})}&=&\sum h_1\smi (\ty ^3g^2)\tqrb
(\tQra )_{<1>}\tx ^3\tz ^2_2\smi ((\tqra (\tQra )_{<0>})_{<1>}\tx
^2\tz ^2_1)\\
&&\ot (\tqra (\tQra )_{<0>})_{<0>}\tx ^1\tz ^1\ty ^1_{<0>}\mfa
_{<0>}\ot h_2\smi (\ty ^2g^1)\tQrb \tz ^3\ty ^1_{<1>}\mfa
_{<1>}\\
{\rm (\ref{rca4},\ref{rca1})}&=&\sum h_1\smi (\ty ^3g^2)\tqrb
\tx ^3(\tQra )_{<1>_2}\smi ((\tqra )_{<1>}\tx ^2(\tQra )_{<1>_1}\b
)\\
&&\ot (\tqra )_{<0>}\tx ^1(\tQra )_{<0>}\ty ^1_{<0>}\mfa _{<0>}\ot
h_2\smi (\ty ^2g^1)\tQrb \ty ^1_{<1>}\mfa _{<1>}\\
{\rm (\ref{rca3},\ref{tpr})}&=&\sum h_1\smi (\ty ^3g^2)\tqrb
\smi ((\tqra )_{<1>}\tprb )\ot (\tqra )_{<0>}\tpra \tQra \ty
^1_{<0>}\mfa _{<0>}\\
&&\ot h_2\smi (\ty ^2g^1)\tQrb \ty ^1_{<1>}\mfa _{<1>}\\
{\rm (\ref{tpqr2},\ref{rca2})}&=&\sum h_1\smi (\ty ^3x ^3\tX
^3_2g^2)\ot \ty ^1\tX ^1\mfa _{<0>}\ot h_2\smi (\a \ty ^2_2x ^2\tX
^3_1g^1)\ty ^2_1x^1\tX ^2\mfa _{<1>}\\
{\rm (\ref{rca3},\ref{ql})}&=&\sum h_1\smi (q^2_L\tX
^3_2g^2)\ot \tX ^1\mfa _{<0>}\ot h_2\smi (q^1_L\tX ^3_1g^1)\tX
^2\mfa _{<1>}\\
&=& \r_{\cal U}(h\ot \mfa )
\end{eqnarray*}
as needed.
\end{examples}

The aim of this Section is to prove that, over a finite
dimensional quasi-Hopf algebra $H$, the category $\hba $ is
isomorphic to a certain category of relative Hopf modules defined
in \cite{bn}. Let $H$ be a finite dimensional
quasi-bialgebra and $A$ a left $H$-module algebra. Recall that a
$k$-vector space $M$ is called a right $(H^*, A)$-Hopf module if
$M$ is a right $H^*$-comodule and a right $A$-module in the
monoidal category
of right $H^*$-comodules ${\cal M}^{H^*}$. In terms of $H$ this means:
\begin{itemize}
\item[-] $M$ is a left $H$-module; denote the action of $h\in H$ on
$m\in M$ by $h\bullet m$;
\item[-] $A$ acts on $M$ from the right; denote the action of $a\in A$ on
$m\in M$ by $m\bullet a$;
\item[-] for all $m\in M$, $h\in H$ and $a, a'\in A$, we have
\begin{eqnarray}
&&m\bullet 1_A=m,\nonumber\\%
&&(m\bullet a)\bullet a'=\sum (X^1\bullet m)\bullet [(X^2\cd
a)(X^3\cd a')],\label{rhm1}\\[1mm]%
&&h\bullet (m\bullet a)=\sum (h_1\bullet m)\bullet (h_2\cd a).
\label{rhm2}
\end{eqnarray}
\end{itemize}
${\cal M}_A^{H^*}$ will be the category of right $(H^*,
A)$-Hopf modules and $A$-linear $H^*$-colinear maps.\\
Before we can establish the claimed isomorphism of categories, we need
a few Lemmas.

\begin{lemma}\lelabel{3.3}
Let $H$ be a finite dimensional quasi-Hopf algebra and $({\mf A},
\r , \Phi _{\r })$ a right $H$-comodule algebra. We have a functor
$$%
F:\ \hba\to {\cal M}_{{\mf A}\ovsm H^*}^{H^*}.%
$$%
For $M\in\hba $, $F(M)=M$, with structure maps
\begin{itemize}
\item[-] $M$ is a left $H$-module via $h\bullet m=S^2(h)\succ m$, $m\in
M$, $h\in H$;
\item[-] ${\mf A}\ovsm H^*$ acts on $M$ from the right by
\end{itemize}
\begin{equation}\label{rhs1}
\hspace{-3cm}m\bullet (\mfa \ovsm \v )=\sum \v \smi
(S(U^1)f^2m_{(1)}\mfa_{<1>}\tprb ) S(U^2)f^1\succ m_{(0)}\prec
\mfa _{<0>}\tpra
\end{equation}
where we denote
\begin{equation}\label{u}
\hspace{-3cm}U=\sum U^1\ot U^2=\sum g^1S(q^2_R)\ot g^2S(q^1_R).
\end{equation}
\end{lemma}

\begin{proof}
The most difficult par of the proof is to show that $F(M)$ satisfies
the relations (\ref{rhm1}) and (\ref{rhm2}). It is then straightforward
to show that a map in $\hba$ is also a map in ${\cal M}_{{\mf A}\ovsm
H^*}^{H^*}$,
and that $F$ is a functor.\\
By \cite[Lemma 3.13]{hn3}
we have, for all $h\in H$:
\begin{eqnarray}
&&U[1\ot S(h)]=\sum \Delta (S(h_1))U(h_2\ot 1),\label{u1}\\%
&&\Phi ^{-1}(id \ot \Delta )(U)(1\ot U)=(\Delta \ot id ) (\Delta
(S(X^1))U)(X^2\ot X^3\ot 1).\label{u2}
\end{eqnarray}
Write $f=\sum f^1\ot f^2=\sum F^1\ot F^2$, $f^{-1}=\sum g^1\ot
g^2$, $\tilde {p}_{\r }=\sum \tpra \ot \tprb =\sum \tPra \ot \tPra
$, and $U=\sum U^1\ot U^2=\sum {\bf U}^1\ot {\bf U}^2$. For all
$m\in M$, $\mfa , \mfa '\in {\mf A}$,
and $\v , \psi \in H^*$, we compute that
\begin{eqnarray*}
&&\hspace*{-2cm}(X^1\bullet m)\bullet \{[X^2\cd (\mfa \ovsm \v )][X^3\cd (\mfa
'\ovsm \psi )]\}\\
&=&\sum (S^2(X^1)\succ m)\bullet [(\mfa \ovsm X^2\rh \v )(\mfa
'\ovsm
X^3 \rh \psi )]\\
&=&\sum (S^2(X^1)\succ m)\bullet [\mfa \mfa '_{<0>}\tx ^1\ovsm
(X^2\rh \v \lh \mfa '_{<1>}\tx ^2)(X^3\rh \psi \lh \tx ^3)]\\
&=&\sum <(X^2\rh \v \lh \mfa '_{<1>}\tx ^2)(X^3\rh \psi \lh \tx
^3), \smi (S(U^1)f^2S^2(X^1)_2m_{(1)}\\
&&(\mfa \mfa '_{<0>}\tx ^1)_{<1>}\tprb )>
S(U^2)f^1S^2(X^1)_1\succ m_{(0)}\prec (\mfa \mfa '_{<0>}\tx
^1)_{<0>}\tpra \\
{\rm (\ref{ca})}&=&\sum <\v , \smi
(F^2S(U^1)_2S(S(X^1)_1)_2f^2_2m_{(1)_2}\mfa _{<1>_2}\mfa '_{<0,
1>_2}\tx ^1_{<1>_2}(\tprb )_2\\
&&g^2S(\mfa '_{<1>}\tx ^2))X^2><\psi , \smi
(F^1S(U^1)_1S(S(X^1)_1)_1f^2_1m_{(1)_1}\mfa _{<1>_1}\\
&&\mfa '_{<0, 1>_1}\tx ^1_{<1>_1}(\tprb )_1g^1S(\tx
^3))X^3>\\
&&S(S(X^1)_2U^2)f^1\succ m_{(0)}\prec
\mfa _{<0>}\mfa '_{<0, 0>}\tx ^1_{<0>}\tpra \\
{\rm (\ref{ca},\ref{tqr2},\ref{rca1})}&=&\sum \v \smi
(S(S(X^1)_{(1, 1)}U^1_1X^2) F^2f^2_2m_{(1)_2}\mfa _{<1>_2}\tX
^3\mfa '_{<0,
1>}\tprb S(\mfa '_{<1>}))\\
&&\psi \smi (S(S(X^1)_{(1, 2)}U^1_2X^3)F^1
f^2_1 m_{(1)_1}\mfa _{<1>_1}\tX^2 %
(\mfa '_{<0, 0>}\tpra )_{<1>}\tPrb )\\
&&S(S(X^1)_2U^2)f^1\succ m_{(0)}\prec \mfa _{<0>}\tX ^1 (\mfa'_{<0, 0>}\tpra
)_{<0>}\tPra \\
{\rm (\ref{u2},\ref{tpqr1})}&=&\sum \v \smi
(S(x^1U^1)F^2f^2_2m_{(1)_2}\mfa _{<1>_2}\tX ^3\tprb )\\
&&\psi \smi (S(x^2U^2_1{\bf U}^1)F^1f^2_1m_{(1)_1}\mfa _{<1>_1}
\tX ^2(\tpra \mfa ')_{<1>}\tPrb )\\
&&S(x^3U^2_2{\bf U}^2)f^1\succ m_{(0)}\prec \mfa _{<0>}\tX
^1(\tpra \mfa ')_{<0>}\tPrb \\
{\rm (\ref{g2},\ref{pf},\ref{rca1})}~~~
&=&\sum \v \smi (S(U^1)F^2m_{(1)}\mfa _{<1>} \tprb )\psi \smi
(S(U^2_1{\bf U}^1)f^2F^1_2m_{(0, 1)}\mfa
_{<0, 1>} \\
&&(\tpra \mfa ')_{<1>}\tPrb )
S(U^2_2{\bf U}^2)f^1F^1_1\succ
m_{(0, 0)}\prec \mfa _{<0, 0>}
(\tpra \mfa ')_{<0>}\tPra \\
{\rm (\ref{ca},\ref{rhs1})}~~~
&=&\sum \v \smi (S(U^1)F^2m_{(1)}
\mfa _{<1>}\tprb )(S(U^2)F^1\succ m_{(0)}%
\prec \mfa _{<0>}\tpra )\bullet (\mfa '\ovsm \psi )\\%
{\rm (\ref{rhs1})}&=&[m\bullet (\mfa \ovsm \v )]\bullet (\mfa
'\ovsm \psi ).
\end{eqnarray*}
Similar computations show that
\begin{eqnarray*}
&&\hspace*{-2cm}\sum (h_1\bullet m)\bullet (h_2\cd (\mfa \ovsm \v ))\\
&=&\sum (S^2(h_1)\succ m)\bullet (\mfa \ovsm h_2\rh \v )\\
&=&\sum <\v , \smi (S(U^1)f^2S^2(h_1)_2m_{(1)}\mfa _{<1>}\tprb
)h_2>\\
&&S(U^2)f^1S^2(h_1)_1\succ m_{(0)}\prec \mfa _{<0>}\tpra \\
{\rm (\ref{ca},\ref{u1})}&=&\sum \v \smi (S(U^1)f^2m_{(1)} \mfa
_{<1>}\tprb )S^2(h)f^1\succ m_{(0)}\prec \mfa _{<0>}\tpra \\
{\rm (\ref{rhs1})}&=&h\bullet [m\bullet (\mfa \ovsm \v )].
\end{eqnarray*}
\end{proof}

Let us next discuss the construction in the converse direction.

\begin{lemma}\lelabel{3.4}
Let $H$ be a finite dimensional quasi-Hopf algebra, $({\mf A}, \r
, \Phi _{\r })$ a right $H$-comodule algebra and $M$ a right
$({\mf A}\ovsm H^*, H^*)$-Hopf module.
Then we have a functor
$$%
G:\ {\cal M}_{{\mf A}\ovsm H^*}^{H^*}\to \hba .%
$$%
For $M\in {\cal M}_{{\mf A}\ovsm H^*}^{H^*}$, $G(M)=M$, with
structure maps
($h\in H$, $m\in M$, $\mfa \in {\mf A}$):
\begin{itemize}
\item[-] $h{\bf \succ}m=S^{-2}(h)\bullet m$;
\item[-] $m{\bf \prec}{\mfa }=m\bullet (\mfa \ovsm \va )$;
\item[-] $\r _M:\ M\ra M\ot H$ is given by
\end{itemize}
\begin{equation}\label{cqhm}
\r_M(m)=\sum m_{\{0\}}\ot m_{\{1\}}=\sum_{i=1}^n [\smi
(V^2g^2)\bullet m]\bullet (\tqra \ovsm \smi (V^1g^1) \rh e^iS\lh
\tqrb )\ot e_i
\end{equation}
where $\{e_i\}_{i=\ov {1, n}}$ and $\{e^i\}_{i=\ov {1, n}}$ are
dual bases and
\begin{equation}\label{v}
V=\sum V^1\ot V^2=\sum \smi (f^2p^2_R)\ot \smi (f^1p^1_R).
\end{equation}
\end{lemma}

\begin{proof}
As in the previous part, the main thing to show is that $G(M)$ is
an object of $\hba$. It is then straightforward to show that $G$
behaves well on the level of the morphisms ($G$ is the identity on
the morphisms).\\%
From the fact that $S^{-2}$ is an algebra map,
it follows that $M$ is a left $H$-module via the action $h\succ
m=S^{-2}(h)\bullet m$. Take the map $$i:\ {\mf A}\ra {\mf A}\ovsm
H^*,~~~i(\mfa )=\mfa \ovsm \va $$ for all $\mfa\in {\mf A}$. Then
$i$ is injective map, $i(1_{\mf A})=1_{{\mf A}\ovsm H^*}$, and
$i(\mfa \mfa ')= i(\mfa )i(\mfa ')$, for all $\mfa , \mfa '\in
{\mf A}$. Therefore, $M$ becomes a right ${\mf A}$-module by
setting $m\prec \mfa =m\bullet i(\mfa )= m\bullet (\mfa \ovsm \va
)$, $m\in M$, $\mfa \in {\mf A}$. Moreover, it is not hard to see
that, with this structure, $M$ is an $(H,{\mf A})$-bimodule. In
order to check the relations (\ref{thm1}-\ref{thm3}) we need some
formulas due to Hausser and Nill \cite[Lemma 3.13]{hn1}, namely
\begin{eqnarray}
&&[1\ot \smi (h)]V=\sum (h_2\ot 1)V\Delta (\smi (h_1)),
\label{v1}\\%
&&(\Delta \ot id)(V)\Phi ^{-1}=\sum (X^2\ot X^3\ot 1)(1\ot V)
(id\ot \Delta )(V\Delta (\smi (X^1))).\label{v2}
\end{eqnarray}
Also, it is clear that
\begin{equation}\label{f1}
(\v \lh h)S=\smi (h)\rh \v S~~~;~~~(h\rh \v )S=\v S\lh
\smi (h)
\end{equation}
for all $h\in H$ and $\v \in H^*$.
Using (\ref{ca}), it follows that
\begin{equation}\label{f2}
(\v S)(\psi S)=\sum [(g^1\rh \psi \lh f^1)(g^2\rh \v \lh f^2)]S
\end{equation}
for all $\v , \psi \in H^*$.
Now, for any $h\in H$ and $m\in M$ we compute that
\begin{eqnarray*}
&&\hspace*{-2cm} \sum h_1\succ m_{\{0\}}\ot h_2m_{\{1\}}\\
&=&\sum_{i=1}^nS^{-2}(h_1)\bullet [(\smi (V^2g^2)\bullet
m)\bullet (\tqra \ovsm \smi (V^1g^1)\rh e^iS\lh \tqrb)]\ot
h_2e_i\\
{\rm (\ref{rhm2})}&=&\sum_{i=1}^n [S^{-2}(h_1)_1\smi
(V^2g^2)\bullet m]\\
&&\bullet (\tqra \ovsm S^{-2}(h_1)_2\smi (V^1g^1)\rh (e^i\lh
h_2)S\lh \tqrb )\ot e_i\\
{\rm (\ref{ca},\ref{f1})}&=&\sum_{i=1}^n [\smi
(V^2\smi (h_1)_2g^2)\bullet m]\\
&&\bullet (\tqra \ovsm \smi (h_2V^1\smi (h_1)_1g^1)
\rh e^iS\lh \tqrb )\ot e_i\\
{\rm (\ref{v1})}&=&\sum_{i=1}^n [\smi
(V^2g^2)S^{-2}(h)\bullet m]\bullet (\tqra \ovsm \smi (V^1g^1)\rh
e^iS\lh \tqrb )\ot e_i\\
&=&\r _M(S^{-2}(h)\bullet m)=\r _M(h\succ m),
\end{eqnarray*}
and similarly, for any $m\in M$ and $\mfa \in {\mf A}$
\begin{eqnarray*}
&&\hspace*{-2cm}
\sum m_{\{0\}}\prec \mfa_{<0>}\ot m_{\{1\}}\mfa _{<1>}\\
&=&\sum_{i=1}^n [(\smi (V^2g^2)\bullet m)\bullet (\tqra
\ovsm \smi (V^1g^1)\rh e^iS\lh \tqrb )]\bullet (\mfa _{<0>}\ovsm
\va )\ot e_i\mfa _{<1>}\\
{\rm (\ref{rhm1},\ref{qsm})}&=&\sum_{i=1}^n (\smi
(V^2g^2)\bullet m)\\
&&\bullet (\tqra \mfa _{<0, 0>}\ovsm \smi (V^1g^1)\rh (\mfa
_{<1>}\rh e^i)S\lh \tqrb \mfa _{<0, 1>})\ot e_i\\
{\rm (\ref{f1},\ref{tpqr1})}&=&\sum_{i=1}^n [\smi
(V^2g^2)\bullet m]\bullet (\mfa \tqra \ovsm \smi (V^1g^1)\rh
e^iS\lh \tqrb )\ot e_i\\
{\rm (\ref{qsm},\ref{rhm1})}&=&\sum_{i=1}^n [(\smi (V^2g^2)\bullet
m)\bullet (\mfa \ovsm \va
)]\bullet (\tqra \ovsm \smi (V^1g^1)\rh e^iS\lh \tqrb )\ot e_i\\
{\rm (\ref{rhm2})}&=&\r _M(m\bullet (\mfa \ovsm \va )) =\r
_M(m\prec \mfa )
\end{eqnarray*}
so the relations (\ref{thm3}) hold. (\ref{thm1}) is obviously
satisfied, thus remain to check (\ref{thm2}) for our structures.
For this we denote $f=\sum f^1\ot f^2$, $f^{-1}=\sum g^1\ot
g^2=\sum G^1\ot G^2$, $\tilde {q}_{\r }=\sum \tqra \ot \tqrb =\sum
\tQra \ot \tQrb $ and $V=\sum V^1\ot V^2=\sum {\bf V}^1\ot {\bf
V}^2$. For all $m\in M$ we compute that
\begin{eqnarray*}
&&\hspace*{-2cm}
X^1\succ m_{\{0, 0\}}\ot X^2m_{\{0, 1\}}\ot
X^3m_{\{1\}}\\
{\rm (\ref{cqhm},\ref{rhm2})}&=&\sum_{i, j=1}^n
\{S^{-2}(X^1)_1\smi ({\bf V}^2G^2)\\
&&\bullet [(\smi (V^2g^2) \bullet m)\bullet (\tqra \ovsm \smi
(V^1g^1)\rh e^iS\lh \tqrb )]\}\\
&&\bullet (\tQra \ovsm S^{-2}(X^1)_2\smi
({\bf V}^1G^1)\rh e^jS\lh \tQrb )\ot X^2e_j\ot X^3e_i\\
{\rm (\ref{ca})}&=&\sum_{i, j=1}^n \{\smi ({\bf V}^2\smi
(X^1)_2G^2)\bullet [(\smi (V^2g^2)\bullet m)\bullet (\tqra \ovsm
\smi (V^1g^1)\\
&&\rh (e^i\lh X^3)S\lh \tqrb )]\} \bullet (\tQra
\ovsm \smi ({\bf V}^1\smi (X^1)_1G^1)\\
&&\rh (e^j\lh X^2)S\lh \tQrb )\ot e_j\ot e_i\\
{\rm (\ref{ca},\ref{rhm2},\ref{f1})}~~~
&=&\sum_{i, j=1}^n\{ [\smi
(V^2{\bf V}^2_2\smi (X^1)_{(2, 2)}G^2_2g^2)\bullet m]\\
&&\bullet (\tqra \ovsm \smi (X^3V^1{\bf V}^2_1\smi (X^1)_{(2,
1)}G^2_1g^1)\rh e^iS\lh \tqrb )\}\\
&&\bullet (\tQra \ovsm \smi (X^2{\bf V}^1 \smi (X^1)_1G^1)\rh
e^jS\lh \tQrb )\ot e_j\ot e_i \\
{\rm (\ref{v2},\ref{rhm1})}&=&\sum_{i, j=1}^n (X^1\smi
(V^2x^3G^2_2g^2)\bullet m)\bullet [(\tqra \ovsm X^2\smi
(V^1_2x^2G^2_1g^1)\rh e^iS \\
&&\rh \tqrb )(\tQra \rh X^3\smi (V^1_1x^1G^1)\rh e^jS\lh \tQrb
)]\ot e_j\ot e_i \\
{\rm (\ref{g2},\ref{pf})}&=&\sum_{i, j=1}^n
(\smi (V^2g^2)\bullet m)\bullet [(\tqra \ovsm \smi
(V^1_2g^1_2G^2)\rh e^iS\lh \tqrb )\\
&&(\tQra \ovsm \smi (V^1_1g^1_1G^1)\rh e^jS\lh \tQrb )]
\ot e_j\ot e_i \\
{\rm (\ref{ca},\ref{qsm})}&=&\sum_{i,j=1}^n (\smi
(V^2g^2)\bullet m)\bullet [\tqra (\tQra )_{<0>}\tx ^1\ovsm (\smi
(V^1g^1)_1\smi (G^2)\rh e^iS\\
&&\lh \tqrb (\tQra )_{<1>}\tx ^2)(\smi (V^1g^1)_2\smi (G^1)\rh
e^jS\lh \tQrb \tx ^3)]\ot e_j\ot e_i\\
{\rm (\ref{tqr2},\ref{mbia2})}&=&\sum_{i, j=1}^n [\smi
(V^2g^2)\bullet m]\bullet \{\tqra \tX ^1_{<0>}\ovsm \smi
(V^1g^1)\rh [(\smi (G^2)\rh e^iS\\
&&\lh \smi (f^2\tX ^3))(\smi (G^1)\rh e^jS\lh \smi (f^1\tX
^2))]\lh \tqrb \tX ^1_{<1>}\} \ot e_j\ot e_i\\
{\rm (\ref{qsm},\ref{f1})}&=&\sum_{i, j=1}^n [\smi
(V^2g^2)\bullet m]\bullet \{[\tqra \ovsm \smi (V^1g^1)\rh (f^2\tX
^3\rh e^i\lh G^2)S\\
&&(f^1\tX ^2\rh e^j\lh G^1)S\lh \tqrb ](\tX
^1\ovsm \va )\} \ot e_j\ot e_i  \\
{\rm (\ref{f2},\ref{rhm1})}&=&\sum_{i, j=1}^n [(\smi
(V^2g^2)\bullet m)\bullet (\tqra \ovsm \smi (V^1g^1)\rh
(e^je^i)S\lh \tqrb )]\\
&&\bullet (\tX ^1\ovsm \va )\ot e_j\tX ^2\ot
e_i\tX ^3\\
{\rm (\ref{cqhm})}&=&\sum m_{\{0\}}\prec \tX ^1 \ot m_{\{1\}_1}\tX
^2\ot m_{\{1\}_2}\tX ^3
\end{eqnarray*}
and this finishes the proof.
\end{proof}

We are now able to prove the main result of this Section,
generalizing \cite[Proposition 2.3]{drvo}.

\begin{theorem}\thlabel{3.5}
Let $H$ be a finite dimensional quasi-Hopf algebra and $({\mf A},
\r , \Phi _{\r })$ a right $H$-comodule algebra. Then the category
of two-sided $(H, {\mf A})$-Hopf modules $\hba $ is isomorphic to
the category of right $(H^*, {\mf A}\ovsm H^*)$-Hopf modules
${\cal M}_{{\mf A}\ovsm H^*}^{H^*}$.
\end{theorem}

\begin{proof}
We have to show that the functors $F$ and $G$ from Lemmas
\ref{le:3.3} and \ref{le:3.4} are inverses.\\%
First, let $M\in \hba $. The structures on $G(F(M))$ (using first
\leref{3.3} and then \leref{3.4}) are denoted by $\succ '$, $\prec
'$ and $\r '_M$. For any $m\in M$, $h\in H$ and $\mfa \in {\mf A}$
we have that
\begin{eqnarray*}
&& h\succ 'm=S^{-2}(h)\bullet
m=S^2(S^{-2}(h))\succ m=h\succ m\\
&&m\prec '\mfa =m\bullet (\mfa \ovsm \va )=m\prec
\mfa
\end{eqnarray*}
because $\sum \va (U^1)U^2=\sum \va (f^2)f^1=1$ and $\sum \va
(m_{(1)})m_{(0)}=m$, $\sum \va (\mfa _{<1>}) \mfa _{<0>}=\mfa $.
In order to prove that $\r'_M=\r _M$, observe first that $\sum
g^1S(g^2\a )=\b $, where we write $f^{-1}=\sum g^1\ot g^2$. This
equality together with (\ref{u},\ref{g2}) and (\ref{pf}) implies
\begin{equation}\label{f3}
\sum g^2_2U^2\ot g^1S(g^2_1U^1)=\sum p^2_L\ot S(p^1_L)
\end{equation}
where $p_L=\sum p^1_L\ot p^2_L$ is the element defined by
(\ref{ql}). Secondly, by $\sum \smi (f^2)\b f^1=\smi (\a )$,
(\ref{g2}) and (\ref{pf}) we have that
\begin{equation}\label{f4}
\sum S(p^2_L)f^1F^1_1\ot \smi (F^2)S(p^1_L)f^2F^1_2=q_R
\end{equation}
were $\sum F^1\ot F^2$ is another copy of $f$, and $q_R$ is the
element defined by (\ref{qr}). Finally, from (\ref{f3},\ref{f4})
and (\ref{pqr}), it follows that
\begin{equation}\label{f5}
\sum S(g^2_2U^2)f^1F^1_1(p_R^1)_1\ot \smi (F^2p^2_R)
g^1S(g^2_1U^1)f^2F^1_2(p^1_R)_2=1\ot 1.
\end{equation}
We now compute for $m\in M$ that
\begin{eqnarray*}
\r '_M(m) &=&\sum_{i=1}^n [\smi (V^2g^2)\bullet
m]\bullet (\tqra \ovsm \smi (V^1g^1)\rh e^iS\lh \tqrb )\ot e_i\\
&=&\sum_{i=1}^n [S(V^2g^2)\succ m]\bullet (\tqra \ovsm
\smi (V^1g^1)\rh e^iS\lh \tqrb )\ot e_i\\
{\rm (\ref{rhs1})}&=&\sum_{i=1}^n <\smi (V^1g^1)\rh
e^iS\lh \tqrb , \smi (S(U^1)f^2S(V^2g^2)_2m_{(1)} (\tqra
)_{<1>}\tprb )>\\
&&S(U^2)f^1S(V^2g^2)_1\succ m_{(0)}\prec (\tqra
)_{<0>}\tpra \ot e_i \\
{\rm (\ref{ca})}&=&\sum S(V^2_2g^2_2U^2)f^1\succ m_{(0)}\prec
(\tqra )_{<0>}\tpra \ot V^1g^1S(V^2_1g^2_1U^1)f^2\\
&&m_{(1)} (\tqra )_{<1>}\tprb S(\tqrb )\\
{\rm (\ref{tpqr2})}&=&\sum S(V^2_2g^2_2U^2)f^1\succ m_{(0)}\ot
V^1g^1S(V^2_1g^2_1U^1)f^2m_{(1)}\\
{\rm (\ref{v},\ref{ca})}&=&\sum S(g^2_2U^2)f^1F^1_1(p^1_R)_1
\succ m_{(0)}\ot \smi (F^2p^2_R)g^1S(g^2_1U^1)f^2
F^1_2(p^1_R)_2m_{(1)} \\
{\rm (\ref{f5})}&=&\sum m_{(0)}\ot m_{(1)}=\r _M(m).
\end{eqnarray*}
and this finishes the proof of the fact that $G(F(M))=M$.\\%

Conversely, take $M\in{\cal M}_{{\mf A}\ovsm
H^*}^{H^*}$. We want to show that $F(G(M))=M$. Denote the left $H$-action
and the right ${\mf A}$-action on $F(G(M))$ by $\bullet'$. Using Lemmas
\ref{le:3.3} and \ref{le:3.4}, we find, for all $h\in H$ and
$m\in M$:%
$$%
h\bullet'm=S^2(h) \succ m=S^{-2}(S^2(h))\bullet m=h\bullet m.%
$$%
The proof of the fact that the right ${\mf A}$-action $\bullet$ and
$\bullet'$ on $M$ coincide is a somewhat more complicated.
Since
$\sum f^2\smi (f^1\b )=S(\a )$, (\ref{g2}) and (\ref{pf}) imply
\begin{equation}\label{f6}
\sum F^1f^1_1p^1_R\ot f^2\smi (F^2f^1_2p^2_R)=\sum S(q^2_L)\ot
q^1_L
\end{equation}
where $q_L=\sum q^1_L\ot q^2_L$ is the element defined by
(\ref{ql}). Also, by (\ref{g2}), (\ref{pf}) and using $\sum S(g^1)\a
g^2=S(\b )$ we can prove the following relation
\begin{equation}\label{f7}
\sum S(G^1)q^1_LG^2_1g^1\ot q^2_LG^2_2g^2=\sum S(p^2_R)\ot
S(p^1_R)
\end{equation}
where $\sum G^1\ot G^2$ is another copy of $f^{-1}$.
Now, from (\ref{u},\ref{ca},\ref{f6},\ref{f7}) and
(\ref{pqr}) it follows that
\begin{equation}\label{f8}
\sum \smi (F^1f^1_1p^1_R)U^2_2g^2\ot S(U^1)f^2\smi
(F^2f^1_2p^2_R)U^2_1g^1=1\ot 1.
\end{equation}
Therefore, for all $m\in M$, $\mfa \in {\mf A}$ and $\v \in H^*$
we have that
\begin{eqnarray*}
&&\hspace*{-2cm}
m\bullet '(\mfa \ovsm \v )\\
{\rm (\ref{rhs1})}&=&\sum \v \smi (S(U^1)f^2m_{\{1\}}\mfa
_{<1>}\tprb ) S(U^2)f^1\succ m_{\{0\}}\prec \mfa _{<0>}\tpra \\
{\rm (\ref{cqhm},\ref{rhm1},\ref{qsm})}&=&
\sum_{i=1}^n \v \smi (S(U^1)f^2e_i\mfa _{<1>}\tprb )S^{-2}(S(U^2)f^1)
\bullet \{[\smi (V^2g^2)\bullet m]\\
&&\bullet [\tqra \mfa _{<0, 0>}(\tpra )_{<0>}\ovsm \smi
(V^1g^1)\rh e^iS\lh \tqrb \mfa _{<0, 1>}(\tpra )_{<1>}]\}\\
&=&\sum_{i=1}^n \v (e_i)S^{-2}(S(U^2)f^1)\bullet \{[\smi
(V^2g^2)\bullet m]\bullet [\tqra \mfa _{<0, 0>}(\tpra )_{<0>}\\
&&\ovsm \smi (V^1g^1)\rh (\mfa _{<1>}\tprb \rh e^i\smi \lh
S(U^1)f^2)S\\
&&\lh \tqrb \mfa _{<0, 1>}(\tpra )_{<1>}]\}\\
{\rm (\ref{f1},\ref{tpqr1},\ref{tpqr2})}&=&\sum
S^{-2}(S(U^2)f^1)\bullet \{[\smi (V^2g^2)\bullet m]\\
&&\bullet [\mfa \ovsm \smi (S(U^1)f^2V^1g^1)\rh \v ]\} \\
{\rm (\ref{rhm2},\ref{ca})}&=&\sum [\smi (V^2\smi
(S(U^2)f^1)_2g^2)\bullet m]\\
&&\bullet [\mfa \ovsm \smi (S(U^1)f^2V^1 \smi
(S(U^2)f^1)_1g^1))\rh \v ]\\
{\rm (\ref{v},\ref{ca})}&=&\sum
[\smi (\smi (F^1f^1_1p^1_R)U^2_2g^2)\bullet m]\\
&&\bullet [\mfa \ovsm \smi (S(U^1)f^2 \smi
(F^2f^1_2p^2_R)U^2_1g^1)\rh \v ]\\
{\rm (\ref{f8})}&=&m\bullet (\mfa \ovsm \v )
\end{eqnarray*}
and this finishes our proof.
\end{proof}

If $H$ is a finite dimensional quasi-Hopf algebra and $A$ is a
left $H$-module algebra then the category ${\cal M}_A^{H^*}$ is
isomorphic to the category of right modules over the smash product
$A\# H$ (\cite[Proposition 2.7]{bn}). Let $M$ be a right $A\# H$-module,
and denote the right action of $a\# h\in A\# H$ on $m\in M$
by $m\rightact (a\# h)$. Following \cite{bn}, $M$ is a right $(H^*,
A)$-Hopf module,
with structure maps

\begin{equation}\label{shm1}
h\bullet m=m\rightact (1\# S(h)), \mbox{${\;\;}$} m\bullet a=\sum
m\rightact [g^1S(q^2_R)\cd a\# g^2S(q^1_R)]
\end{equation}
for all $m\in M$, $a\in A$ and $h\in H$. Conversely, if $M$ is a
right $(H^*, A)$-Hopf module then $M$ is a right $A\#
H$-module, with $A\# H$-action
\begin{equation}\label{ssmm}
m\rightact (a\# h)=\sum \smi (h)\bullet [(\smi (q^2_Lg^2)\bullet
m)\bullet (\smi (q^1_Lg^1)\cd a)].
\end{equation}
Here $q_R=\sum q^1_R\ot q^2_R$, $q_L=\sum q^1_L\ot q^2_L$ and
$f^{-1}=\sum g^1\ot g^2$ are the elements defined by (\ref{qr}),
(\ref{ql}) and (\ref{g}). Combining this with \thref{3.5}, we
obtain the following result.

\begin{corollary}\colabel{3.6}
Let $H$ be a finite dimensional quasi-Hopf algebra and $({\mf A},
\r , \Phi _{\r })$ a right $H$-comodule algebra. Then the category
$\hba $ is isomorphic to the category of right $({\mf A}\ovsm
H^*)\# H$-modules, ${\cal M}_{({\mf A}\ovsm H^*)\# H}$.
\end{corollary}

For later use, we describe the isomorphism of \coref{3.6} explicitely,
leaving verification of the details to the reader.\\
First take $M\in {\cal M}_{({\mf A}\ovsm H^*)\#H}$. The following structure
maps make $M\in \hba$:
\begin{eqnarray}
&&h\succ m=m\rightact ((1_{\mf A}\ovsm \va)\# \smi (h))\label{sth1}\\
&&m\prec \mfa =m\rightact ((\mfa \ovsm \va )\# 1)\label{sth1a}\\
&&\r _M(m)=\sum_{i=1}^n [(\tqra \ovsm \smi (g^2)\rh e^iS\lh
\tqrb )\# \smi (g^1)]\ot e_i\label{sth2}
\end{eqnarray}
for all $m\in M$, $h\in H$ and $\mfa \in {\mf A}$.
$\tilde {q}_{\r }=\sum \tqra \ot
\tqrb $ is the element defined in (\ref{tqr}), $\{e_i\}$ is a basis of $H$ and
$\{e^i\}$ is the corresponding dual basis of $H^*$.\\
Now take $M\in \hba $. Then $M$ is a right $({\mf A}\ovsm H^*)\#H$-module
via the action
\begin{equation}\label{ss2}
m\rightact [(\mfa \ovsm \v )\# h]=\sum \v \smi (f^2m_{(1)}\mfa
_{<1>} \tprb )S(h)f^1\succ m_{(0)}\prec \mfa_{<0>}\tpra .
\end{equation}
In \cite{hn3}, it is shown that, for a finite dimensional
quasi-Hopf algebra $H$, the
category of right quasi-Hopf $H$-bimodules ${}_H{\cal M}_H^H$
naturally coincides with the category of representations of the
two-sided crossed product $H\gtl H^*\trl H$ constructed in
\cite{hn1}. We will show in \seref{4} that $H\gtl H^*\trl
H=(H\ovsm H^*)\# H$ as algebras.

%%%%%%%%%%%%%%%%%%%%%%%%%%%%%%%%%%%%%%%%%%%%%%%%%%%%%%%%
\section{Two-sided crossed products are generalized smash products}\selabel{4}
%%%%%%%%%%%%%%%%%%%%%%%%%%%%%%%%%%%%%%%%%%%%%%%%%%%%%%%%
\setcounter{equation}{0}
Let $H$ be a finite dimensional quasi-bialgebra, and $(\mf {A}, \r ,
\Phi _{\r })$,$(\mf {B}, \l ,
\Phi _{\l })$ respectively a right and a left $H$-comodule algebra.
As in the case of a Hopf algebra,
the right $H$-coaction $(\r , \Phi _{\r })$ on $\mf{A}$ induces
a left $H^*$-action
$\tr :\ H^*\ot \mf {A}\ra \mf {A}$ given by
\begin{equation}\label{tra}
\varphi \tr \mf {a}=\sum \varphi (\mf {a}_{<1>})\mf {a}_{<0>}
\end{equation}
for all $\varphi \in H^*$ and $\mf {a}\in \mf{A}$, and
where $\r (\mf {a})=\sum a_{<0>}\ot \mf {a}_{<1>}$ for
any $\mf {a}\in \mf {A}$. Similarly, the left $H$-action $(\l ,
\Phi _{\l })$ on $\mf {B}$ provides a right $H^*$-action
$\tl :
\mf{B}\ot H^*\ra \mf {B}$ given by
\begin{equation}\label{tla}
\mf {b}\tl \varphi =\sum \varphi (\mf {b}_{[-1]})\mf {b}_{[0]}
\end{equation}
for all $\varphi \in H^*$ and $\mf {b}\in \mf {B}$, where we now
denote $\l (\mf {b})=\sum \mf {b}_{[-1]}\ot \mf {b}_{[0]}$ for
$\mf {b}\in \mf {B}$. Following \cite[Proposition 11.4 (ii)]{hn1}
we can define an algebra structure on the $k$-vector space $\mf
{A}\ot H^*\ot \mf {B}$. This algebra is denoted by $\mf {A}\gtl
_{\r }H^*\trl _{\l }\mf {B}$ and its multiplication is given by
\begin{eqnarray}
&&\hspace*{-2cm}
(\mf {a}\gtl \varphi \trl
\mf {b})(\mf {a}'\gtl \psi \trl \mf {b}')\nonumber\\
&=&\sum \mf {a}(\varphi _1\tr \mf {a}')\tx ^1_{\r }\gtl (\tx
^1_{\l }\rh \varphi _2\lh \tx ^2_{\r })(\tx ^2_{\l }\rh \psi
_1\lh \tx ^3_{\r })\trl \tx ^3_{\l }(\mf {b}\tl \psi _2)\mf {b}'\label{tscp}
\end{eqnarray}
for all $\mf {a}, \mf {a}'\in \mf {A}$, $\mf {b}, \mf
{b}'\in \mf {B}$, and $\v, \psi \in H^*$,
where we write $\mf {a}\gtl \varphi \trl
\mf {b}$ for $\mf {a}\ot \varphi \ot \mf {b}$
when viewed as an element of $\mf {A}\gtl _{\r }H^*\trl _{\l }\mf {B}$. The
unit of
the algebra $\mf {A}\gtl_{\r }H^*\trl_{\l }\mf {B}$ is $1_{\mf
{A}}\gtl \va \trl 1_{\mf {B}}$. Hausser and Nill called this algebra the
two-sided crossed product.
In this Section we will prove that this two-sided crossed product algebra
is a generalized smash product
between the quasi-smash product ${\mf A}\ovsm H^*$ and ${\mf B}$.

\begin{proposition}\prlabel{4.1}
Let $H$ be a quasi-bialgebra, $A$ a left $H$-module algebra and
$\mathfrak{B}$ a left $H$-comodule algebra. Let
$A\gsm \mathfrak{B}=A\ot\mathfrak{B}$ as a $k$-module, with newly defined
multiplication
\begin{equation}\label{gsm}
(a\gsm \mf {b})(a'\gsm \mf {b}')=\sum (\tilde {x}^1\cd
a)(\tilde {x}^2\mf {b}_{[-1]}\cd a')\gsm \tilde {x}^3\mf {
b}_{[0]}\mf {b}'
\end{equation}
for all $a, a'\in A$ and $\mf {b}, \mf {b}'\in \mathfrak{B}$.
Then $A\gsm \mathfrak{B}$ is an associative algebra
with unit $1_A\gsm 1_{\mathfrak{B}}$.
\end{proposition}

\begin{proof}
For all $a, a', a''\in A$ and $\mf {b}, \mf {b}', %
\mf {b}''\in \mathfrak{B}$ we have:
\begin{eqnarray*}
&&\hspace*{-2cm}
[(a\gsm b)(a'\gsm \mf {b}')](a''\gsm
\mf {b}'')\\
&=&\sum [(\tilde {x}^1\cd a)(\tilde {x}^2\mf {b}_{[-1]}\cd
a')\gsm \tilde {x}^3\mf {b}_{[0]}
\mf {b}'](a''\gsm \mf {b}'')\\
&=&\sum [(\ty ^1_1\tx ^1\cd a)(\ty ^1_2\tx ^2\mf {b}_{[-1]}\cd a')]
(\ty ^2\tx ^3_{[-1]}\mf {b}_{[0, -1]}\mf {b}'_{[-1]}\cd a'')\gsm
\ty ^3\tx ^3_{[0]}\mf {b}_{[0, 0]}\mf {b}'_{[0]}\mf {b}''\\
{\rm (\ref{ma1})}&=&\sum (X^1\ty ^1_1\tx ^1\cd a)[(X^2\ty
^1_2\tx ^2\mf {b}_{[-1]}\cd a')(X^3\ty ^2\tx ^3_{[-1]}
\mf {b}_{[0, -1]}\mf {b}'_{[-1]}\cd a'')]\\
&&\gsm \ty ^3\tx ^3_{[0]}\mf {b}_{[0, 0]}\mf {b}'_{[0]}\mf {b}''\\
{\rm (\ref{lca2})}&=&\sum (\tx ^1\cd a)[(\tx ^2_1\ty ^1\mf
{b}_{[-1]}\cd a')(\tx ^2_2\ty ^2\mf {b}_{[0, -1]}\mf
{b}'_{[-1]}\cd a'')]\gsm \tx ^3\ty ^3\mf {b}_{[0,
0]}\mf {b}'_{[0]}\mf {b}''\\
{\rm (\ref{lca1},\ref{ma2})}&=&\sum (\tx ^1\cd a)\{(\tx
^2\mf {b}_{[-1]}\cd [(\ty ^1\cd
a')(\ty ^2\mf {b}'_{[-1]}\cd a'')]\}\gsm
\tx ^3\mf {b}_{[0]}\ty ^3\mf {b}'_{[0]}\mf {b}''\\
&=&\sum (a\gsm \mf {b})[(\ty ^1\cd a')(\ty ^2\mf
{b}'_{[-1]}\cd a'')\gsm \ty ^3\mf {b}'_{[0]}\mf
{b}'']\\
&=&(a\gsm \mf {b})[(a'\gsm \mf {b}')(a''\gsm \mf
{b}'')].
\end{eqnarray*}
It follows from
(\ref{lca3}), (\ref{lca4}) and (\ref{ma3}) that
$1_A\gsm 1_{\mf {B}}$ is the unit for $A\gsm \mf {B}$.
\end{proof}

\begin{remark}\relabel{4.2}
Let $H$ be a quasi-bialgebra and $A$ a left $H$-module algebra.
Then $H$ is a left $H$-comodule algebra so it make sense to
consider $A\gsm H$. It is not hard to see that in this case $A\gsm
H$ is just the smash product $A\# H$. For this
reason we will call the algebra $A\gsm \mf {B}$ in \prref{4.1}
the generalized smash product of $A$ and $\mf {B}$. In fact,
our terminology is in agreement with the terminology used over Hopf algebras,
see \cite{doi} and \cite{cmz}.
\end{remark}

Let $H$ be a finite dimensional quasi-bialgebra, $(\mf
{A}, \r , \Phi _{\r })$ a right $H$-comodule algebra and $(\mf
{B}, \l , \Phi _{\l })$ a left $H$-comodule algebra. Then the
quasi-smash product $\mf {A}\ovsm H^*$ is a left $H$-module
algebra so it makes sense to consider the generalized smash
product $(\mf {A}\ovsm H^*)\gsm \mf {B}$. The main result of this
Section is now the following:

\begin{proposition}\prlabel{4.3}
With notation as above, the algebras $(\mf {A}\ovsm H^*)\gsm
\mf {B} $ and $\mf {A}\gtl _{\r }H^*\trl _{\l }\mf {B}$ coincide.
\end{proposition}

\begin{proof}
Using (\ref{gsm}), (\ref{aqsm}) and (\ref{qsm}) we compute that the
multiplication on $(\mf {A}\ov {\#}H^*)\gsm \mf {B}$
is given by
\begin{eqnarray*}
&&\hspace*{-2cm}[(\mfa \ovsm \v )\gsm \mfb ][(\mfa '\ovsm \psi )\gsm \mfb']\\
&=&\sum [\tx ^1_{\l }\cd (\mfa \ovsm \v )][\tx ^2_{\l }\mfb
_{[-1]}\cd (\mfa '\ovsm \psi )]\gsm \tx ^3_{\l }\mfb _{[0]}\mfb
'\\
&=&\sum (\mfa \ovsm \tx ^1_{\l }\rh \v )(\mfa '\ovsm \tx ^2_{\l
}\mfb _{[-1]}\rh \psi )\gsm \tx ^3_{\l }\mfb _{[0]}\mfb '\\
&=&\sum \mfa \mfa '_{<0>}\tx ^1_{\r }\ovsm (\tx ^1_{\l }\rh \v
\lh \mfa '_{<1>}\tx ^2_{\r })(\tx ^2_{\l }\mfb _{[-1]}\rh \psi
\lh \tx ^3_{\r })\gsm \tx ^3_{\l }\mfb _{[0]}\mfb '\\
{\rm (\ref{tra},\ref{tla})}&=&\sum \mfa (\v _1\tr \mfa
')\tx ^1_{\r }\ovsm (\tx ^1_{\l }\rh \v _2\lh \tx ^2_{\r })
(\tx ^2_{\l }\rh \psi _1\lh \tx ^3_{\r })\gsm \tx ^3_{\l }(\mfb
'\tl \psi _2)\mfb'
\end{eqnarray*}
for $\mfa , \mfa '\in \mf {A}$, $\mfb , \mfb '\in \mf
{B}$, and $\v , \psi \in H^*$. This is just the multiplication
rule on the two-sided crossed product $\mf {A}\gtl_{\r }H^*\trl
_{\l }\mf {B}$.
\end{proof}

Take $\mf {B}=H$ in \prref{4.3}. From \reref{4.2}, we obtain:

\begin{corollary}\colabel{4.4}
Let $H$ be a finite dimensional quasi-bialgebra and $(\mf {A}, \r
, \Phi _{\r })$ a right $H$-comodule algebra. Then $({\mf A}\ovsm
H^*)\# H=\mf {A}\gtl _{\r }H^*\trl _{\Delta }H$ as algebras. In
particular, $(H\ovsm H^*)\# H=H\gtl H^*\trl H$ as algebras.
\end{corollary}

%%%%%%%%%%%%%%%%%%%%%%%%%%%%%%%%%%%%%%%%%%%%%%%%%%%%%%%%%%%%%%%
\section{The category of Doi-Hopf modules}\selabel{5}
%%%%%%%%%%%%%%%%%%%%%%%%%%%%%%%%%%%%%%%%%%%%%%%%%%%%%%%%%%%%%%%
\setcounter{equation}{0}
Let $H$ be a Hopf algebra over a field $k$, $A$ an $H$-comodule
algebra and $C$ an $H$-module coalgebra. A Doi-Hopf module
is a $k$-vector space together with an $A$-action and a $C$-coaction
satisfying a certain compatibility relation. They were introduced
independently by Doi \cite{doi} and Koppinen \cite{Koppinen95}, and it turns
out that most types of Hopf modules that had been studied before
were special cases: Sweedler's Hopf modules
\cite{sw}, Doi's relative Hopf modules \cite{doi2}, Takeuchi's
relative Hopf modules \cite{ta}, Yetter-Drinfeld modules,
graded modules and  modules graded by a $G$-set.\\
Over a quasi-Hopf algebra, the category of relative Hopf modules has been
introduced
and studied \cite{bn}, as well as the category of Hopf $H$-bimodules (see
\cite{hn3})
and the category of Hopf modules $_H^H{\cal M}_H^H$ (see \cite{sh}).
We will introduce Doi-Hopf modules, and
we will show that, at least in the case where $H$ is finite dimensional,
all these
categories are isomorphic to certain categories of Doi-Hopf
modules.\\%
First we recall from \cite{bn} the definition
of a relative Hopf module.
Let $H$ be a quasi-bialgebra and $C$ a right $H$-module coalgebra.
Let $N$ be a $k$-vector space furnished with the following additional
structure:
\begin{itemize}
\item[-] $N$ is a right $H$-module; the
right action of $h\in H$ on $n\in N$ is denoted by $nh$;
\item[-] $N$ is a left $C$-comodule in the
monoidal category ${\cal M}_H$; we use the following notation for the
left $C$-coaction on $N$: $\r_N:\ N\ra C\ot N$, $\r _N(n)=\sum
n_{[-1]}\ot n_{[0]}$; this means that the following conditions hold,
for all $n\in N$:
\begin{eqnarray}
&&\sum \une (n_{[-1]})n_{[0]}=n\nonumber\\
&&(\und \ot id_N)(\r _N(n))\Phi ^{-1}=(id_C\ot \r _N)(\r _N(n))\label{lhm1}
\end{eqnarray}
\item[-] we have the following compatibility relation, for all $n\in N$ and
$c\in C$:
\begin{equation}
\r _N(nh)=\sum n_{[-1]}\cd h_1\ot n_{[0]}h_2.\label{lhm2}
\end{equation}
\end{itemize}
Then $N$ is called a left $[C,H]$-Hopf module. $^C{\cal M}_H$ is the
category of left
$[C, H]$-Hopf modules; the morphisms are right $H$-linear maps which are
also left $C$-comodule maps. We will now generalize this definition.

\begin{definition}\delabel{5.1}
Let $H$ be a quasi-bialgebra over a field $k$, $C$ a right
$H$-module coalgebra and $(\mf {B}, \l , \Phi _{\l })$ a left
$H$-comodule algebra. A right-left $(H, \mf {B}, C)$-Hopf module
(or Doi-Hopf module) is a $k$-module $N$, with the following additional
structure: $N$ is
right $\mf {B}$-module (the right action of
$\mfb $ on $n$ is denoted by $n\mfb $), and we have a $k$-linear map
$\r_N:\ N\ra C\ot N$, such that the following relations hold, for
all $n\in N$ and $\mfb \in \mf {B}$:
\begin{eqnarray}
&&(\und \ot id_N)(\r _N(n))=(id_C\ot \r _N) (\r _N(n))\Phi _{\l }
\label{dhm1}\\%
&&(\une \ot id_N)(\r _N(n))=n\label{dhm2}\\%
&&\r _N(n\mfb )=\sum n_{[-1]}\cd \mfb _{[-1]}\ot n_{[0]}\mfb
_{[0]}.\label{dhm3}
\end{eqnarray}
As usual, we use the Sweedler-type notation $\r _N(n)= \sum n_{[-1]}\ot
n_{[0]}$.
$^C{\cal M}(H)_{\mf {B}}$ is the category of right-left $(H,
{\mf B}, C)$-Hopf modules and right ${\mf B}$-linear, left
$C$-colinear $k$-linear maps.
\end{definition}

Obviously, if $\mf {B}=H$, $\l =\Delta $ and $\Phi _{\l }=\Phi $,
then $^C{\cal M}(H)_{\mf {B}}= {}^C{\cal M}_H$.\\
The main aim of \seref{6} will be to define the category
of two-sided two-cosided Hopf modules over a quasi-bialgebra, and to
prove that it is isomorphic to a module category in the finite dimensional
case. To this end, we will need our next result, stating
that the category of Doi-Hopf modules is
a module category in the case where the coalgebra $C$ is finite dimensional.
In fact, for an
arbitrary right $H$-module coalgebra $C$, the linear dual space of
$C$, $C^*$, is a left $H$-module algebra. The multiplication of
$C^*$ is the convolution, that is $(c^*d^*)(c) =\sum c^*(\una
)d^*(\unb )$, the unit is $\une $ and the left $H$-module
structure is given by $(h\rh c^*)(c)=c^*(c\cd h)$, for $h\in H$, $c^*,
d^*\in C^*$, $c\in C$. Thus $C^*$ is a left $H$-module algebra and
$(\mf {B}, \l , \Phi _{\l })$ is a left $H$-comodule algebra. By
\prref{4.1}, it makes sense to consider the generalized smash
product algebra $C^*\gsm \mf {B}$.

\begin{proposition}\prlabel{5.2}
Let $H$ be a quasi-bialgebra, $C$ a finite dimensional right
$H$-module coalgebra and $(\mf {B}, \l , \Phi _{\l})$ a left
$H$-comodule algebra. Then the category $^C{\cal M}(H)_{\mf {B}}$
of right-left $(H, \mf {B}, C)$-Hopf modules is isomorphic to the
category ${\cal M}_{C^*\gsm {\mf B}}$
of right modules over $C^*\gsm {\mf B}$.
\end{proposition}

\begin{proof}
We restrict to defining the functors that define the isomorphism
of categories, leaving all other details to the reader.
Let $\{c_i\}_{i=\ov{1, n}}$ and
$\{c^i\}_{i=\ov{1, n}}$ be dual bases in $C$ and $C^*$.\\
Let $N$ be a right $C^*\gsm \mf {B}$-module. Since $i:\ \mf {B}\ra
C^*\gsm \mf {B}$, $i(\mfb )=\une \gsm \mfb $ for $\mfb \in \mf
{B}$, is an algebra map, it follows that $N$ is a right $\mf
{B}$-module via the action $n\mfb =ni(\mfb )=n(\une \gsm \mfb )$,
$n\in N$, $\mfb \in \mf {B}$. The map $j:\ C^*\ra C^*\gsm \mf
{B}$, $j(c^*)=c^*\gsm 1_{\mf {B}}$, $c^*\in C^*$, is not an
algebra map (it is not multiplicative) but it can be used to define a left
$C$-coaction on $N$:
\begin{equation}\label{c1}
\r _N(n)=\sum n_{[-1]}\ot n_{[-1]}= \sum_{i=1}^nc_i\ot
nj(c^*)=\sum_{i=1}^nc_i\ot n(c^i\gsm 1_{\mf {B}}).
\end{equation}
We can easily check that $N$ becomes an
object in ${}^C{\cal M}(H)_{\mf {B}}$.\\
Conversely, take $N\in {}^C{\cal M}(H)_{\mf {B}}$. Then $N$ is a right
$\mf {B}$-module and $C^*$ acts on $M$ from the right as follows:
let $nc^*=\sum
c^*(n_{-1]})n_{[0]}$, $n\in N$, $c^*\in C^*$. Now define
\begin{equation}\label{c2}
n(c^*\gsm \mfb )=(nc^*)\mfb =\sum c^*(n_{[-1]})n_{[0]}\mfb .
\end{equation}
Then $N$ becomes a right $C^*\gsm \mf
{B}$-module.
\end{proof}

%%%%%%%%%%%%%%%%%%%%%%%%%%%%%%%%%%%%%%%%%%%%%%%%%%%
\section{Two-sided two-cosided Hopf modules}\selabel{6}
%%%%%%%%%%%%%%%%%%%%%%%%%%%%%%%%%%%%%%%%%%%%%%%%%%%
\setcounter{equation}{0}
Now we define the category of two-sided two-cosided
Hopf modules ${}_H^C{\cal M}_{\mathbb {A}}^H$. If
$H$ is finite dimensional, then this category is isomorphic to a certain
category of right-left Doi-Hopf modules, $^C{\cal M}(H\ot
H^{\rm op})_{(\mb {A}\ovsm H^*)\# H}$. As a consequence, if $C$ is
also finite dimensional then this category is isomorphic to the
category of right modules over a generalized smash product, by
\prref{5.2}.

\begin{definition}\delabel{6.1}\cite[Definition 8.2]{hn1}.
Let $H$ be a quasi-bialgebra. An $H$-bicomodule algebra $\mb
{A}$ is a quintuple $(\mb{A},\l, \r , \Phi_{\l }, \Phi_{\r }, \Phi
_{\l , \r })$, where $\l $ and $\r $ are left and right
$H$-coactions on $\mb {A}$,
 and where $\Phi _{\l}\in H\ot H\ot \mb {A}$, $\Phi _{\r }\in \mb {A}\ot
H\ot H$ and
$\Phi _{\l , \r }\in H\ot \mb {A}\ot H$ are invertible elements,
such that
\begin{itemize}
\item[-] $(\mb {A}, \l , \Phi _{\l })$ is a left $H$-comodule algebra,
\item[-] $(\mb {A}, \r , \Phi _{\r })$ is a right $H$-comodule algebra,
\item[-] the following compatibility relations hold, for all $a\in \mb
{A}$:
\end{itemize}${\;\;}$\\[-1.5cm]
\begin{eqnarray}
&&\hspace{-1cm}\Phi _{\l , \r }(\l \ot id)(\r (a))=(id\ot \r )(\l
(a))\Phi _{\l, \r }\label{bca1}\\%
&&\hspace{-1cm}(1_H\ot \Phi _{\l
, \r })(id\ot \l \ot id)(\Phi _{\l , \r }) (\Phi _{\l }\ot
1_H)=(id\ot id\ot \r )(\Phi _{\l })(\Delta \ot id\ot id)(\Phi _{\l
, \r }) \label{bca2}\\%
&&\hspace{-1cm}(1_H\ot \Phi _{\r })(id\ot
\r \ot id)(\Phi _{\l ,\r })(\Phi _{\l , \r }\ot 1_H)= (id\ot id\ot
\Delta )(\Phi _{\l , \r })(\l \ot id\ot id) (\Phi _{\r
}).\label{bca3}
\end{eqnarray}
\end{definition}

It was pointed out in \cite{hn1} that the following additional
relations hold in an $H$-bicomodule algebra
$\mb {A}$:
\begin{equation}\label{bca4}
(id_H\ot id_{\mb {A}}\ot \va )(\Phi _{\l , \r })=1_H\ot 1_{\mb
{A}}, \mbox{${\;\;}$} (\va \ot id_{\mb {A}}\ot id_H)(\Phi _{\l ,
\r })= 1_{\mb {A}} \ot 1_H.
\end{equation}
As a first example, take $\mb {A}=H$, $\l =\r =\Delta $ and $\Phi
_{\l }=\Phi _{\r }= \Phi _{\l , \r }=\Phi $. Related to the left
and right comodule algebra structures of $\mb {A}$ we will keep
the notation of the previous Sections. We will use  the following
notation: %
$$%
\Phi _{\l , \r }=\sum \O ^1\ot \O ^2\ot \O ^3=\sum
\ov {\O }^1\ot \ov {\O }^2\ot \ov {\O }^3={\rm etc.} %
$$%
and %
$$%
\Phi ^{-1}_{\l , \r }=\sum \o ^1\ot \o ^2\ot \o ^3=\sum \ov {\o
}^1\ot \ov {\o }^2\ot \ov {\o }^3={\rm etc.} %
$$%
If $H$ is a quasi-bialgebra, then the opposite algebra $H^{\rm
op}$ is also a quasi-bialgebra. The reassociator of $H^{\rm op}$
is $\Phi _{\rm op}=\Phi ^{-1}$. $H\ot H^{\rm op}$ is also a
quasi-bialgebra with reassociator
\begin{equation}\label{phhop}
\Phi_{H\ot H^{\rm op}}=\sum (X^1\ot x^1)\ot (X^2\ot x^2)\ot
(X^3\ot x^3).
\end{equation}
If we identify $H\ot H^{\rm op}$-modules and $(H,H)$-bimodules, then
the category of $(H,H)$-bimodules, ${}_H{\cal M}_H$, is monoidal. The
associativity constraints are given by ${\bf a'}_{U, V, W}:
(U\ot V)\ot W\ra U\ot (V\ot W)$, where
\begin{equation}\label{bim}
{\bf a'}_{U, V, W}((u\ot v)\ot w)= \Phi \cd (u\ot (v\ot w))\cd\Phi ^{-1}
\end{equation}
for all $U, V, W\in {}_H{\cal M}_H$, $u\in U$, $v\in V$ and $w\in
W$. A coalgebra in the category of
$(H,H)$-bimodules will be called an $H$-bimodule
coalgebra. More precisely, an $H$-bimodule coalgebra $C$ is an
$(H,H)$-bimodule (denote the actions by $h\cd c$ and $c\cd h$)
with a comultiplication $\und :\ C\ra C\ot C$ and
a counit $\une :\ C\ra k$ satisfying the following relations,
for all $c\in C$ and $h\in H$:
\begin{eqnarray}
&&\Phi \cd (\und \ot id)(\und (c))\cd \Phi
^{-1}=(id\ot \und
)(\und (c))\label{bmc1}\\
&&\und (h\cd c)=\sum h_1\cd \una \ot h_2\cd \unb ,
\mbox{${\;\;}$} \und (c\cd h)=\sum \una \cd h_1\ot \unb \cd h_2
\label{bmc2}\\
&&\une (h\cd c)=\va (h)\une (c), \mbox{${\;\;}$} \une
(c\cd h)=\une (c)\va (h)\label{bmc3}
\end{eqnarray}
where we used the same Sweedler-type notation as before.
An $H$-bimodule coalgebra $C$ becomes a right $H\ot H^{\rm op}$-module
coalgebra via
the right $H\ot H^{\rm op}$-action
\begin{equation}\label{rmcs}
c\cd (h\ot h')=h'\cd c\cd h
\end{equation}
for $c\in C$ and $h, h'\in H$.
Our next definition extends the definition of two-sided two-cosided Hopf
modules
from \cite{sh}.

\begin{definition}\delabel{6.2}
Let $H$ be a quasi-bialgebra, $(\mb {A}, \l , \r , \Phi _{\l },
\Phi _{\r }, \Phi _{\l , \r })$ an $H$-bicomodule algebra and $C$ an
$H$-bimodule coalgebra. A two-sided
two-cosided $(H, \mb {A}, C)$-Hopf module (or crossed Hopf module)
is a $k$-vector space with the following additional structure:
\begin{itemize}
\item[-] $N$ is an $(H, \mb {A})$-two-sided Hopf module, i.e. $N\in
{}_H{\cal M}_{\mb {A}}^H$; we write $\succ $ for the left $H$-action,
$\prec $ for the right ${\mb A}$-action, and
$\r_N^H(n)= \sum n_{(0)}\ot n_{(1)}$ for the right
$H$-coaction on $n\in N$;
\item[-] we have $k$-linear map $\r _N^C:\ N\ra C\ot N$,
$\r _N^C(n)=\sum n_{[-1]}\ot n_{[0]}$, called the left $C$-coaction
on $N$, such that $\sum \une (n_{[-1]})n_{[0]}=n$ and
\begin{equation}\label{tstc1}
\Phi (\und \ot id_N)(\r ^C_N(n))=(id_C\ot \r ^C_N) (\r
^C_N(n))\Phi _{\l }
\end{equation}
for all $n\in N$;
\item[-] $N$ is a $(C,H)$-``bicomodule", in the sense that, for all $n\in N$,
\begin{equation}\label{tstc2}
\Phi (\r ^C_N\ot id_H)(\r _N^H(n))=(id_C\ot \r ^H_N) (\r
^C_N(n))\Phi _{\l , \r }
\end{equation}
\item[-] the following compatibility relations hold
\begin{eqnarray}
\r ^C_N(h\succ n)&=&\sum h_1\cd n_{[-1]}\ot h_2\succ n_{[0]}\label{tstc3}\\
\r ^H_N(n)(n\prec a)&=&\sum n_{[-1]}\cd a_{[-1]}\ot
n_{[0]}\prec a_{[0]}\label{tstc3b}
\end{eqnarray}
for all $h\in H$, $n\in N$ and $a\in \mb {A}$.
\end{itemize}
${}^C_H{\cal M}_{\mb {A}}^H$ will be the category of two-sided
two-cosided Hopf modules and maps preserving the actions by $H$
and $\mb {A}$ and the coactions by $H$ and $C$.
\end{definition}

Let $H$ be a quasi-bialgebra, $\mb {A}$ an $H$-bicomodule algebra
and $C$ an $H$-bimodule coalgebra. Let us call the threetuple $(H, \mb {A},
C)$ a
{\sl Drinfeld datum}. In the rest of this Section we will show that if
$H$ is a finite dimensional quasi-Hopf algebra then the above
category ${}_H^C{\cal M}_{\mb {A}}^H$ is isomorphic to a certain
category of Doi-Hopf modules. To this end, we first need some
lemmas.

\begin{lemma}\lelabel{6.3}
Let $H$ be a finite dimensional quasi-Hopf algebra and $\ba $ an
$H$-bicomodule algebra. Consider the map
$$\wp :\ (\mbA \ovsm H^*)\#
H\ra (H\ot H^{\rm op})\ot (\mbA \ovsm H^*)\# H$$
 given by
\begin{equation}\label{lcas1}
\wp ((a\ovsm \v )\# h)=\sum a_{[-1]}\o ^1\ot S(y^3h_2)\ot
(a_{[0]}\o ^2\ovsm y^1\rh \v \lh \o ^3)\# y^2h_1
\end{equation}
for any $a\in \mbA $, $\v \in H^*$ and $h\in H$, where $\Phi
^{-1}_{\l , \r }=\sum \o ^1\ot \o ^2\ot \o ^3$. Set
\begin{equation}\label{lcas2}
\Phi _{\wp }=\sum (\tX _{\l }^1\ot g^1S(x^3))\ot (\tX ^2_{\l }\ot
g^2S(x^2))\ot (\tX ^3_{\l }\ovsm \va )\# x^1
\end{equation}
where $f^{-1}=\sum g^1\ot g^2$ is the element defined in
(\ref{g}). Then $((\mbA \ovsm H^*)\# H, \wp , \Phi _{\wp})$ is a
left $H\ot H^{\rm op}$-comodule algebra.
\end{lemma}

\begin{proof}
We first show that $\wp $ is an algebra map. Using
(\ref{sm1}) and (\ref{qsm}) we can easily show that the
multiplication on $(\mb A \ovsm H^*)\# H$ is given by
\begin{eqnarray}
&&\hspace*{-2cm}((a\ovsm \v )\# h)((a'\ovsm \psi )\# h')\nonumber\\
&=& \sum
aa'_{<0>}\tx ^1_{\r }\ovsm (x^1\rh \v \lh a'_{<1>}\tx ^2_{\r
})(x^2h_1\rh \psi \lh \tx ^3_{\r }) \# x^3h_2h'\label{smdp}
\end{eqnarray}
for all $a, a'\in \mb A$, $\v , \psi \in H^*$ and $h, h'\in
H$. Therefore
\begin{eqnarray*}
&&\hspace*{-2cm}\wp (((a\ovsm \v )\# h)((a'\ovsm \psi )\# h'))\\
&=&\sum a_{[-1]}a'_{<0>_{[-1]}}(\tx ^1_{\r })_{[-1]}\o ^1 \ot
S(y^3x^3_2h_{(2, 2)}h'_2)\ot [a_{[0]}a'_{<0>_{[0]}}(\tx
^1_{\r })_{[0]}\o ^2\\
&&\ovsm (y^1_1x^1\rh \v \lh a'_{<1>}\tx ^2_{\r }\o ^3_1)
(y^1_2x^2h_1\rh \psi \lh \tx ^3_{\r }\o ^3_2)]\# y^2x^3_1h_{(2,
1)}h'_1\\
{\rm (\ref{bca3},\ref{q3})}&=&\sum
a_{[-1]}a'_{<0>_{[-1]}}\ov {\o }^1\o ^1\ot S(y^3x^3h_{(2,
2)}h'_2)\ot [a_{[0]}a'_{<0>_{[0]}}\ov {\o }^2\o ^2_{<0>}\tx
^1_{\r }\\
&&\ovsm (z^1y^1\rh \v \lh a'_{<1>}\ov {\o }^3\o ^2_{<1>} \tx
^2_{\r })(z^2y^2_1x^1h_1\rh \psi \lh \o ^3\tx ^3_{\r })]\#
z^3y^2_2x^2h_{(2, 1)}h'_1\\
{\rm (\ref{bca1},\ref{q1})}&=&\sum a_{[-1]}\ov {\o
}^1a'_{[-1]}\o ^1\ot S(y^3h_2)\cdot _{op}S(x^3h'_2)\ot
[a_{[0]}\ov {\o }^2 (a'_{[0]}\o ^2)_{<0>}\tx ^1_{\r }\ovsm
(z^1y^1\\
&&\rh \v \lh \ov {\o }^3(a'_{[0]}\o ^2)_{<1>} \tx ^2_{\r
})(z^2y^2_1h_{(1, 1)}x^1\rh \psi \lh \o ^3\tx ^3_{\r })]\#
z^3y^2_2h_{(1, 2)}x^2h'_1\\
{\rm (\ref{qsm})}&=& \sum a_{[-1]}\ov {\o }^1a'_{[-1]}\o ^1\ot
S(y^3h_2)\cdot _{op}S(x^3h'_2)\ot [(a_{[0]}\ov {\o }^2\ovsm
z^1y^1\rh \v \lh \ov {\o }^3)\\
&&(a'_{[0]}\o ^2\ovsm z^2y^2_1 h_{(1, 1)}x^1\rh \psi \lh \o
^3)]\# z^3y^2_2 h_{(1, 2)}x^2h'_1\\
{\rm (\ref{sm1})}&=& \sum a_{[-1]}\ov {\o }^1a'_{[-1]}\o ^1\ot
S(y^3h_2)\cd _{op}S(x^3h'_2)\ot [(a_{[0]}\ov {\o }^2\ovsm
y^1\rh \v \lh \ov {\o }^3)\# y^2h_1]\\
&&[(a'_{[0]}\o ^2\ovsm x^1\rh \psi \lh \o ^3)\# x^2h'_1]\\
&=&\wp ((a\ovsm \v)\# h)\wp ((a'\ovsm \psi )\# h')
\end{eqnarray*}
where $\cdot_{\rm op}$ is the product in $H^{\rm op}$. Obviously
$\wp $ respects the unit element and (\ref{lca3}) holds.
(\ref{lca1}) can be proved using similar
computations as above and is left to the reader. Using the notation
$$\Phi _{\wp }=\sum \tX ^1_{\wp }\ot \tX ^2_{\wp }\ot \tX
^3_{\wp }={\rm etc.} $$
we can compute
\begin{eqnarray*}
&&\hspace*{-2cm}(id\ot id
\ot \wp )(\Phi _{\wp })(\Delta \ot id\ot id) (\Phi {\wp })\\
&=&\sum (\tX ^1_{\l }\ot g^1S(x^3))((\tY ^1_{\l })_1\ot
G^1_1S(y^3)_1)\ot (\tX ^2_{\l }\ot g^2S(x^2))\\
&&((\tY ^1_{\l })_2\ot G^1_2S(y^3)_2)\ot ((\tX ^3_{\l })_{[-1]}\ot
S(x^1_2))(\tY ^2_{\l }\ot G^2S(y^2))\\
&&\ot [((\tX ^3_{\l })_{[0]}\ovsm \va )\# x^1_1][(\tY ^3_{\l
}\ovsm \va )\# y^1]\\
{\rm (\ref{ca},\ref{q3})}&=& \sum (\tX ^1_{\l }(\tY ^1_{\l
})_1\ot G^1_1g^1S(y^3x^3))\ot (\tX ^2_{\l }(\tY ^1_{\l })_2\ot
G^1_2g^2S(z^3y^2_2x^2))\\
&&\ot ((\tX ^3_{\l })_{[-1]}\tY ^2_{\l }\ot G^2S(z^2y^2_1x^1))\ot
[((\tX ^3_{\l })_{[0]}\tY ^3_{\l }\ovsm \va )\# z^1y^1]\\
{\rm (\ref{lca2},\ref{g2},\ref{pf})}&=&\sum (\tY ^1_{\l
}X^1\ot x^1g^1S(y^3))\ot (\tX ^1_{\l } (\tY ^2_{\l })_1X^2\ot
x^2g^2_1G^1S(z^3y^2_2))\\
&&\ot (\tX ^2_{\l }(\tY ^2_{\l })_2X^3\ot
x^3g^2_2G^2S(z^2y^2_1))\ot
[(\tX ^3_{\l }\tY ^3_{\l }\ovsm \va )\# z^1y^1]\\
{\rm (\ref{ca})}&=&\sum (\tY ^1_{\l }\ot g^1S(y^3))(X^1\ot x^1)\ot
(\tX ^1_{\l }\ot G^1S(z^3))((\tY ^2_{\l })_1 \ot
g^2_1S(y^2)_1)\\
&&(X^2\ot x^2)\ot (\tX ^2_{\l }\ot G^2S(z^2))((\tY ^2_{\l })_2\ot
g^2_2S(y^2)_2)(X^3\ot x^3)\\
&&\ot [(\tX ^3_{\l }\ovsm \va )\# z^1][(\tY ^3_{\l }
\ovsm \va )\# y^1]\\
{\rm (\ref{phhop})}&=&(1_H\ot \Phi _{\wp }) (id\ot \Delta _{H\ot
H^{\rm op}}\ot id)(\Phi _{\wp })(\Phi _{H\ot H^{\rm op}} \ot {\bf 1})
\end{eqnarray*}
where $\sum G^1\ot G^2$ is another copy of $f^{-1}$ and
${\bf 1}=(1_{\mbA }\ovsm \va )\# 1_H$ is the unit of the algebra
$(\mbA \ovsm H^*)\# H$.
\end{proof}

Let $H$ be a finite dimensional quasi-Hopf
algebra, $\ba $ an $H$-bicomodule algebra and $C$ an $H$-bimodule
coalgebra. By \leref{6.3}, we can consider the category of Doi-Hopf
modules ${}^C{\cal M}(H\ot H^{\rm op})_{(\mbA \ovsm H^*)\# H}$.
We will prove that it is isomorphic
to the category of two-sided two-cosided Hopf modules $\tstc $.

\begin{lemma}\lelabel{6.4}
Let $H$ be a quasi-Hopf algebra and $(H, \mbA ,C)$ a Drinfeld datum.
We have a functor
$$%
F:\ \tstc\to {}^C{\cal M}(H\ot H^{\rm op})_{(\mbA \ovsm H^*)\#
H}.%
$$%
$F(N)=N$ as a $k$-module, with structure maps given by the
equations
\begin{eqnarray}
&&n\rightact ((a\ovsm \v )\# h)=\sum  \v \smi
(f^2n_{(1)}a_{<1>}\tprb )S(h)f^1\succ n_{(0)}\prec a_{<0>}\tpra
\label{s1dhm}\\
&&{\tilde {\r }}_N^C(n)=\sum n_{\{-1\}}\ot n_{\{0\}}=\sum f^1\cd
n_{[-1]}\ot f^2\succ n_{[0]}\label{s2dhm}
\end{eqnarray}
for all $n\in N$, $a\in \mbA $, $\v \in H^*$ and $h\in H$.
$F$ sends a morphism to itself.
\end{lemma}

\begin{proof}
Since $N$ is a two-sided $({\mb A}, H)$-Hopf module, we know by
(\ref{ss2}) that $N$ is a right $(\mbA \ovsm H^*)\# H$-module via
the action defined by (\ref{s1dhm}). Let $\sum F^1\ot
F^2$ be another copy of $f$. For any $n\in N$, we have that
\begin{eqnarray*}
&&\hspace*{-2cm}(\un {\Delta }\ot id_N)(\tilde {\r }_N^C(n))\Phi ^{-1}_{\wp }\\
{\rm (\ref{lcas2})}&=&\sum n_{{\{-1\}}_{\un {1}}}\cd (\tx ^1_{\l
}\ot S(X^3)F^1)\ot n_{{\{-1\}}_{\un {2}}}\cd (\tx ^2_{\l }\ot
S(X^2)F^2)\\
&&\ot n_{\{0\}} \rightact [(\tx ^3_{\l }\ovsm \va )\# X^1]\\
{\rm (\ref{rmcs},\ref{s2dhm})}&=&\sum S(X^3)F^1\cd (f^1\cd
n_{[-1]})_{\un {1}}\cd \tx ^1_{\l }\ot S(X^2)F^2\cd (f^1\cd
n_{[-1]})_{\un {2}}\cd \tx ^2_{\l }\\
&&\ot S(X^1)f^2\succ n_{[0]}\prec \tx ^3_{\l }\\
{\rm (\ref{bmc2})}&=&\sum S(X^3)F^1f^1_1\cd n_{[-1]_{\un {1}}}\cd
\tx ^1_{\l }\ot S(X^2)F^2f^1_2\cd n_{[-1]_{\un {2}}}\cd \tx ^2_{\l
}\\
&&\ot S(X^1)f^2\succ n_{[0]}\prec \tx ^3_{\l }\\
{\rm (\ref{tstc1},\ref{g2},\ref{pf})}&=&\sum f^1\cd
n_{[-1]}\ot F^1f^2_1\cd n_{[0, -1]}\ot F^2f^2_2\succ
n_{[0, 0]}\\
{\rm (\ref{tstc3})}&=&\sum f^1\cd n_{[-1]}\ot F^1\cd (f^2\succ
n_{[0]})_{[-1]}\ot F^2\succ (f^2\succ n_{[0]})_{[0]}\\
{\rm (\ref{s2dhm})}&=&\sum n_{\{-1\}}\ot F^1\cd
n_{\{0\}_{[-1]}}\ot F^2\succ n_{\{0\}_{[0]}}\\%
{\rm (\ref{s2dhm})}&=&(id_C\ot \tilde {\r }_N^C)(\tilde {\r
}_N^C(n)).
\end{eqnarray*}
We still have to show the compatibility relation (\ref{dhm3}). First
observe that (\ref{tpr}), (\ref{bca3}) and (\ref{q5}) imply
\begin{equation}\label{fx}
\sum \O ^1(\tpra )_{[-1]}\ot \O ^2(\tpra )_{[0]}\ot \O ^3\tprb
=\sum \o ^1\ot \o ^2_{<0>}\tpra \ot \o ^2_{<1>}\tprb S(\o ^3).
\end{equation}
For all $n\in N$, $a\in \mbA $, $\v \in H^*$ and $h\in H$ we
compute that
\begin{eqnarray*}
&&\hspace*{-2cm}\tilde {\r }_N^C(n\rightact ((a\ovsm \v )\# h))\\
{\rm (\ref{s1dhm},\ref{s2dhm})}&=&\sum \v \smi
(f^2n_{(1)}a_{<1>}\tprb )F^1\cd (S(h)f^1\succ
n_{(0)}\prec a_{<0>}\tpra )_{[-1]}\\
&&\ot F^2\succ (S(h)f^1\succ n_{(0)}\prec a_{<0>}\tpra )_{[0]}\\
{\rm (\ref{tstc3},\ref{ca})}&=&\sum \v \smi
(f^2n_{(1)}a_{<1>}\tprb )S(h_2)F^1f^1_1\cd
n_{(0)_{[-1]}}\cd a_{<0>_{[-1]}}(\tpra )_{[-1]}\\
&&\ot S(h_1)F^2f^1_2\succ n_{(0)_{[0]}}\prec a_{<0>_{[0]}}(\tpra
)_{[0]}\\
{\rm (\ref{bca1})}&=&\sum \v \smi (f^2n_{(1)}\o ^3a_{[0]_{<1>}}\O
^3\tprb )S(h_2)F^1f^1_1\cd n_{(0)_{[-1]}}\cd \o ^1a_{[-1]}\O
^1(\tpra )_{[-1]}\\
&&\ot S(h_1)F^2f^1_2\succ n_{(0)_{[0]}}\prec \o ^2a_{[0]_{<0>}}\O
^2(\tpra )_{[0]}\\
{\rm (\ref{fx},\ref{tstc2})}&=&\sum \v \smi
(f^2x^3n_{[0]_{(1)}}a_{[0]_{<1>}}\o ^2_{<1>}\tprb S(\o
^3))S(h_2)F^1f^1_1x^1\cd n_{[-1]}\cd a_{[-1]}\o ^1\\
&&\ot S(h_1)F^2f^1_2x^2\succ n_{[0]_{(0)}}\prec a_{[0]_{<0>}}\o
^2_{<0>}\tpra \\
{\rm (\ref{g2},\ref{pf})}&=&\sum (x^1\rh \v \lh \o ^3)\smi
(F^2f^2_2n_{[0]_{(1)}}a_{[0]_{<1>}}\o ^2_{<1>}\tprb
)S(x^3h_2)f^1\cd n_{[-1]}\cd a_{[-1]}\o ^1\\
&&\ot S(x^2h_1)F^1f^2_1\succ n_{[0]_{(0)}}\prec a_{[0]_{<0>}}\o
^2_{<0>}\tpra \\
{\rm (\ref{tstc3})}&=&\sum (x^1\rh \v \lh \o ^3)\smi (F^2(f^2\succ
n_{[0]})_{(1)}a_{[0]_{<1>}}\o ^2_{<1>}\tprb )\\
&&S(x^3h_2)f^1\cd n_{[-1]}\cd a_{[-1]}\o ^1\ot S(x^2h_1)F^1\succ
(f^2\succ n_{[0]})_{(0)}\prec a_{[0]_{<0>}}\o ^2_{<0>}\tpra \\
{\rm (\ref{rmcs},\ref{s1dhm})}&=&\sum f^1\cd n_{[-1]}\cd
[a_{[-1]}\o ^1\ot S(x^3h_2)]\\
&&\ot (f^2\succ n_{[0]})\rightact [(a_{[0]}\o ^2\ovsm x^1\rh \v
\lh \o ^3)\# x^2h_1]\\
{\rm (\ref{s2dhm},\ref{lcas1})}&=&\tilde {\r }_N^C(n)\wp
((a\ovsm \v )\# h)
\end{eqnarray*}
completing the proof.
\end{proof}

\begin{lemma}\lelabel{6.5}
Let $H$ be a finite dimensional quasi-Hopf algebra and $(H, \mbA ,
C)$ a Drinfeld datum. We have a functor
$$%
G:\ {}^C{\cal M}(H\ot H^{\rm op})_{({\mbA}\ovsm H^*)\# H}\to
{}^C_H{\cal M}_{\mbA }^H.%
$$%
$G(N)=N$ as a $k$-module, with structure maps given by
\begin{eqnarray}
&&\hspace{-7mm}h\succ n=n\rightact [(1_{\mbA }\ovsm \va )\# \smi
(h)],
\mbox{${\;\;}$}
n\prec a=n\rightact [(a\ovsm \va )\# 1_H],\label{ststc1}\\
&&\hspace{-7mm}\r _N^H:\ N\ra N\ot H, \mbox{${\;\;}$}
\r _N^H(n)=\sum_{i=1}^nn\rightact [(\tqra \ovsm \smi
(g^2)\rh e^iS\lh \tqrb )\# \smi (g^1)]\ot
e_i,\label{ststc2}\\%
&&\hspace{-7mm}\un {\r }_N^C: N\ra C\ot N, \mbox{${\;\;}$}
\un {\r }_N^C(n)=\sum g^1\cd n_{[-1]}\ot g^2\succ
n_{[0]}\label{ststc3}
\end{eqnarray}
for $n\in N$, $a\in \mbA$ and $h\in H$. Here $\{e_i\}_{i=\ov {1,
n}}$ is a basis of $H$ and $\{e^i\}_{i=\ov {1, n}}$ is the
corresponding dual basis of $H^*$. $G$ sends a morphism to itself.
\end{lemma}

\begin{proof}
Since $N$ is a right $(\mbA \ovsm H^*)\# H$-module, we already
know by (\ref{sth1}) and (\ref{sth2}) that $H$ is a two-sided
$(\mbA , H)$-Hopf module via (\ref{ststc1}) and (\ref{ststc2}).
Thus we only have to check (\ref{tstc1}), (\ref{tstc2}) and
(\ref{tstc3}). First note that $N\in {}^C{\cal
M}(H\ot H^{\rm op})_{(\mbA \ovsm H^*)\# H}$ implies
\begin{eqnarray}
&&\hspace*{-2cm} \sum n_{[-1]}\ot n_{[0, -1]}\ot n_{[0,
0]}\nonumber\\%
&=& \sum S(X^3)f^1\cd n_{[-1]_{\un {1}}}\cd \tx
^1_{\l }\ot S(X^2)f^2\cd n_{[-1]_{\un {2}}}\cd \tx ^2_{\l }\ot
n_{[0]}\rightact [(\tx ^3_{\l }\ovsm \va )\# X^1]\label{fx2}\\%
&&\hspace*{-2cm} \sum \{n\rightact [(a\ovsm \v )\# h]\}_{[-1]}\ot
\{n\rightact [(a\ovsm \va )\# h]\}_{[0]}\nonumber\\%
&=& \sum S(x^3h_2)\cd n_{[-1]}\cd a_{[-1]}\o ^1\ot
n_{[0]}\rightact [(a_{[0]}\o ^2\ovsm x^1\rh \v \lh \o ^3)\#
x^2h_1]\label{fx3}
\end{eqnarray}
for all $n\in N$, $a\in \mbA $, $\v \in H^*$ and $h\in H$. By the
above definitions and (\ref{fx3}) it is immediate that
\begin{equation}\label{fx4}
{\un \r }^C_N(h\succ n)=\Delta (h){\un \r }_N^C(n)~~{\rm and}~~
{\un \r }_N^C(n\prec a)={\un \r }_N^C(n)\r _{\l }(a)
\end{equation}
for all $h\in H$, $n\in N$ and $a\in \mbA $
(we leave it to the reader to verify the details). Let
$\sum G^1\ot G^2$ be another copy of $f^{-1}$. We
compute that
\begin{eqnarray*}
&&\hspace*{-2cm}\Phi ({\un \Delta }\ot id_N)(\r _N^C(n))\\
{\rm (\ref{ststc3})}&=&\sum X^1\cd (g^1\cd n_{[-1]})_{\un {1}}\ot
X^2\cd (g^1\cd n_{[-1]})_{\un {2}}\ot X^3g^2\succ n_{[0]}\\
{\rm (\ref{ststc1},\ref{bmc2})}&=&\sum X^1g^1_1\cd n_{[-1]_{\un
1}}\ot X^2g^1_2\cd n_{[-1]_{\un 2}}\ot
n_{[0]}\rightact [(1_{\mbA }\ovsm \va )\# \smi (X^3g^2)]\\
{\rm (\ref{fx2},\ref{smdp})}&=&\sum X^1g^1_1G^1S(x^3)\cd
n_{[-1]}\cd \tX ^1_{\l }\ot X^2g^1_2G^2S(x^2)\cd n_{[0, -1]}\cd
\tX ^2_{\l }\\
&&\ot n_{[0, 0]}\rightact [(\tX ^3_{\l }\ovsm \va )\# \smi
(X^3g^2S(x^1))]\\
{\rm (\ref{g2},\ref{pf})}&=&\sum g^1\cd n_{[-1]}\cd \tX ^1_{\l
}\ot g^2_1G^1\cd n_{[0, -1]}\cd \tX ^2_{\l }\\
&&\ot n_{[0, 0]}\rightact [(\tX ^3_{\l }\ovsm \va )
\# \smi (g^2_2G^2)]\\
{\rm (\ref{ststc1})}&=&\sum g^1\cd n_{[-1]}\cd \tX ^1_{\l }\ot
g^2_1G^1\cd n_{[0, -1]}\cd \tX ^2_{\l }\ot g^2_2G^2\succ n_{[0,
0]}\prec \tX ^3_{\l }\\%
{\rm (\ref{ststc3},\ref{bmc2},\ref{fx4})}&=&(id_C\ot {\un \r
}_N^C)({\un \r }_N^C(n))\Phi _{\l }.
\end{eqnarray*}
The verification of (\ref{tstc2}) is based on
similar computations, and we leave the details to the reader.
\end{proof}

As a consequence of Lemmas \ref{le:6.4} and \ref{le:6.5}, we have
the following description of $^C_H{\cal
M}_{\mbA }^H$ as a category of Doi-Hopf modules;
this description generalizes \cite[Proposition 2.3]{bdr}.

\begin{theorem}\thlabel{6.6}
Let $H$ be a finite dimensional quasi-Hopf algebra and $(H, \mbA ,
C)$ a Drinfeld datum. Then the categories $\tstc $ and $^C{\cal
M}(H\ot H^{\rm op})_{(\mbA \ovsm H^*)\# H}$ are isomorphic.
\end{theorem}

\begin{proof}
We have to verify that the functors $F$ and $G$
defined in Lemmas \ref{le:6.4} and \ref{le:6.5} are inverses.
For the $C$-coactions  (\ref{s2dhm}) and
(\ref{ststc3}), this is obvious; for the other structures, it has
been already done in \coref{3.6}.
\end{proof}

\prref{5.2} and \thref{6.6} immediately imply the following result.

\begin{corollary}\colabel{6.7}
Let $H$ be a finite dimensional quasi-Hopf algebra and $(H, \mbA ,
C)$ a Drinfeld datum with $C$ finite dimensional. Then the
category $\tstc $ is isomorphic to the category of right modules
over the generalized smash product $C^*\gsm ((\mbA \ovsm H^*)\#
H)$.
\end{corollary}

\begin{remark}\relabel{6.8}
Let $H$ be a finite dimensional Hopf algebra. Cibils and Rosso
\cite{cr} introduced an algebra $X=(H^{\rm op}\ot H){\un \ot }(H^*\ot
H^{*op})$ having the property that the category of two-sided two-cosided Hopf
modules over $H^*$ coincides with the category of left $X$-modules.
Moreover, it was also proved in \cite{cr} that $X$ is isomorphic
to the direct tensor product of a Heisenberg double and the
opposite of a Drinfeld double. Recently, Panaite \cite{p}
introduced two other algebras $Y$ and $Z$ with the same property
as $X$. More precisely, $Y$ is the two-sided crossed product $H^*\#
(H\ot H^{\rm op})\# H^{*\rm op}$, and $Z$ is the diagonal crossed product
in the sense of \cite{hn1}, $(H^*\ot H^{*op})\Join (H\ot H^{\rm op})$.
Using different methods, we proved that
the category of two-sided two-cosided Hopf modules over a finite dimensional
quasi-Hopf algebra is isomorphic to the
category of right (resp. left) modules over the generalized smash
product ${\cal A}=H^*\gsm ((H\ovsm H^*)\# H)$ (resp. ${\cal A}^{\rm op}$).
Note that, in general, the multiplication on
$C^*\gsm ((\mbA \ovsm H^*)\# H)$ is given by the formula
\begin{eqnarray*}
&&\hspace*{-1.5cm}[c^*\gsm ((a\ovsm \v )\# h)][d^*\gsm ((a'\ovsm
\psi)\#
h')]\nonumber\\%
&=&\sum (\tx ^1_{\l }\rh c^*\lh S(X^3)f^1) (\tx ^2_{\l }a_{[-1]}\o
^1 \rh d^*\lh S(X^2x^3h_2)f^2) \gsm \{[\tx ^3_{\l }a_{[0]}\o
^2a'_{<0>}\tx ^1_{\r }\nonumber\\
&&
\ovsm (X^1_{(1, 1)}y^1x^1\rh \v \lh \o ^3a'_{<1>}\tx ^2_{\r })
(X^1_{(1, 2)}y^2x^2_1h_{(1, 1)}\rh \psi \lh \tx ^3_{\r })]\#
X^1_2y^3x^2_2h_{(1, 2)}h'\}\label{fx5}
\end{eqnarray*}
\end{remark}
%%%%%%%%%%%%%%%%%%%%%%%%%%%%%%%%%%%%%%%%%%%%%%%%%%%%%%%%%%%%

\end{document}